\documentclass[12pt,leqno]{amsart}

\usepackage[dvips]{graphicx}
\usepackage{amsfonts}
\usepackage{amssymb}
\usepackage{amsmath}
\usepackage{amscd}

\def\DATE{\today}
\newtheorem{theorem}{Theorem}
\newtheorem{definition}[theorem]{Definition}
\newtheorem{corollary}[theorem]{Corollary}

\newtheorem{lemma}[theorem]{Lemma}

\newtheorem{proposition}[theorem]{Proposition}

\usepackage{graphicx}
\email{nicolas.goze@uha.fr, elisabeth.remm@uha.fr}

\begin{document}
\title{An algebraic approach to the set of intervals.}
\author{Nicolas Goze, Elisabeth Remm}
\address{Universit\'e de Haute Alsace, LMIA, 4 rue des Fr\`eres Lumi\`ere,
68093 Mulhouse} \maketitle

\begin{abstract}
In this paper we present the set of intervals as a normed vector space. We define also a four-dimensional associative
algebra whose product gives the product of intervals in any cases. This approach allows to give a notion of divisibility
and in some cases an euclidian division.
\end{abstract}

\section{\protect\bigskip Intervals and generalized intervals}

\bigskip

An interval is a connected closed subset of $\mathbb{R}.$ The classical arithmetic operations on intervals are defined
such that the result of the corresponding operation on elements belonging to operand intervals belongs to the resulting
interval. That is, if $\diamond $ denotes one of the classical operation $+,-,\ast $, we have
\begin{eqnarray}
\lbrack x^{-},x^{+}]\diamond \lbrack y^{-},y^{+}]=\{x\diamond y\text{ }/ \text{ }x\in \lbrack x^{-},x^{+}],\text{ }y\in
\lbrack y^{-},y^{+}]\}\text{ \ \ } .  \label{op}
\end{eqnarray}
In particular we have
\begin{equation*}
\left\{
\begin{array}{c}
\lbrack x^{-},x^{+}]+[y^{-},y^{+}]=[x^{-}+y^{-},x^{+}+y^{+}], \\ \lbrack
x^{-},x^{+}]-[y^{-},y^{+}]=[x^{-}-y^{+},x^{+}-y^{-}]
\end{array}
\right.
\end{equation*}
and
\begin{equation*}
\lbrack x^{-},x^{+}]-[x^{-},x^{+}]=[x^{-}-x^{+},x^{+}-x^{-}]\neq 0.
\end{equation*}
Let $\mathbb{IR}$ be the set of intervals. It is in one to one correspondence with the half plane of $\mathbb{R}^{2}$:
\begin{equation*}
\mathcal{P}_{1}=\{(a,b),a\leq b\}.
\end{equation*}
This set is closed for the addition and $\mathcal{P}_{1}$ is endowed with a regular
semigroup structure. Let $\mathcal{P}_{2}$ be the half plane symmetric to $\mathcal{P}_{1}$ with respect to the first
bisector $\Delta $ of equation $y-x=0.$ The substraction on $\mathbb{IR}$, which is not the symmetric operation of $+$,
corresponds to the following operation on $ \mathcal{P}_{1}$:
\begin{equation*}
(a,b)-(c,d)=(a,b)+s_{\Delta }\circ s_{0}(c,d),
\end{equation*}
where $s_{0}$ is the symmetry with respect to $0$, and $s_{\Delta }$ with respect to $\Delta .$ The multiplication $\ast
$ is not globally defined. Consider the following subset of $\mathcal{P}_{1}$:
\begin{equation*}
\left\{
\begin{array}{l}
\mathcal{P}_{1,1}=\{(a,b) \in \mathcal{P}_{1},a\geq 0, b\geq 0 \}, \\ \mathcal{P}_{1,2}=\{(a,b) \in \mathcal{P}_{1},
a\leq 0, b\geq 0 \}, \\ \mathcal{P}_{1,3}=\{(a,b) \in \mathcal{P}_{1}, a\leq 0, b\leq 0 \}. \\
\end{array}
\right.
\end{equation*}
We have the following cases:

1) If $(a,b),(c,d)\in \mathcal{P}_{1,1}$ the product is written $(a,b)\ast (c,d)=(ac,bd).$

\noindent

\noindent Then if $\ e_{1}=(1,1)$ and $e_{2}=(0,1)$, these "vectors" generate $\mathcal{P}_{1,1}:$
\begin{equation*}
\forall (x,y)\in \mathcal{P}_{1,1} \text{ then }(x,y)=xe_{1}+(y-x)e_{2}, \text{ }x>0,\text{ }y-x>0.
\end{equation*}
The multiplication corresponds in this case to the following associative commutative algebra:
\begin{equation*}
\left\{
\begin{array}{l}
e_{1}e_{1}=e_{1}, \\ e_{1}e_{2}=e_{2}e_{1}=e_{2}e_{2}=e_{2}.
\end{array}
\right.
\end{equation*}

2) Assume that $(a,b)\in \mathcal{P}_{1,1}$ and $(c,d)\in \mathcal{P}_{1,2}$ so $c\leq 0$ and $d\geq 0.$ Thus we obtain
$(a,b)\ast (c,d)=(bc,bd)$ and this product does not depend of $a.$ Then we obtain the same result for any $%
a<b$. Let
$e_{1}=(0,1)$ and $e_{2}=(-1,0).$ Any interval of $\mathcal{P} _{1,1}$ is written $ae_{1}+be_{2}$ with $b<0$ and any
interval of $\mathcal{P }_{1,2}$, $ce_{1}+de_{2}$ with $c,d>0. $ We have no associative multiplication between
$(e_{1},e_{2})$ which describes the product. We have to add a formal dimension to obtain a $3$-dimensional associative
algebra and the product appears as the projection in the plane $(e_{1},e_{2})$ of this associative algebra. Here if we
consider the following associative commutative algebra
\begin{equation*}
\left\{
\begin{array}{l}
e_{1}e_{1}=e_{1},\text{ }e_{1}e_{2}=e_{2},\text{ }e_{1}e_{3}=-e_{2}, \\ e_{2}e_{3}=-e_{1},e_{2}e_{2}=e_{1}, \\
e_{3}e_{3}=e_{3}.
\end{array}
\right.
\end{equation*}
then $(\alpha e_{1}+\beta e_{2}+\beta e_{3})(\gamma e_{1}+\delta e_{2})=\alpha \gamma e_{1}+\alpha \delta e_{2}$. As
$(a,b)=be_{1}-ae_{2}$ and $(c,d)=de_{1}-ce_{2},$ we obtain the expected product.

3) If $(a,b)\in \mathcal{P}_{1,1}$ and $(c,d)\in \mathcal{P}_{1,3}$ then $ a\geq 0,b\geq 0$ and $c\leq 0,d\leq0$ and we
have $(a,b)\ast (c,d)=(bd,ac).$ Let $e_{1}=(1,1)$, $e_{2}=(0,1).$ This product corresponds to the following associative
algebra:
\begin{equation*}
\left\{
\begin{array}{l}
e_{1}e_{1}=e_{1}, \\ e_{1}e_{2}=e_{1}-e_{2}, \\ e_{2}e_{2}=e_{1}-e_{2}.
\end{array}
\right.
\end{equation*}
We have similar results for the cases $(\mathcal{P}_{1,2},\mathcal{P} _{1,2}),\ (\mathcal{P}_{1,2},\mathcal{P}_{1,3})$
and $(\mathcal{P}_{1,3}, \mathcal{P}_{1,3}).$

All this shows that the set $\mathbb{IR}$ is not algebraically structured. Let us describe a vectorial structure on
$\mathbb{IR}$ using the previous geometrical interpretation of $\mathbb{IR}$ with $\mathcal{P}_{1}$. First we extend
$\mathcal{P}_{1}$ to $\mathbb{R}^{2}$ and we obtain an extended set $ \overline{\mathbb{IR}}$ which corresponds to the
classical interval $[a,b]$ and "generalized intervals" $[a,b]$ with $a>b.$ Of course using the addition of
$\mathbb{R}^{2},$ we obtain on $\overline{\mathbb{IR}}$ a structure of abelian group and the symmetric of $[a,b]\in
\mathbb{IR}$ is $[-a,-b]\in \overline{\mathbb{IR}} \setminus \mathbb{IR}.$ In this context $ [a,b]+[-a,-b]=0.$ This
aspect as been developed in \cite{Ma}.

We have a group homomorphism $\varphi $ on $\overline{\mathbb{IR}}$ given by
\begin{equation*}
\begin{tabular}{llll}
$\varphi :$ & $\overline{\mathbb{IR}}$ & $\longrightarrow $ & $\overline{ \mathbb{IR}}$ \\ & $(a,b)$ & $\longrightarrow
$ & $(b,a).$
\end{tabular}
\end{equation*}
This map is called dual and we denote by dual $(a,b)$ the generalized interval $(b,a).$ The corresponding arithmetic has
been developed by Kaucher \cite{Ka} and is naturally called the Kaucher arithmetic. In the following we recall how to
complete the semigroup $\mathbb{IR}$ to obtain a natural vectorial structure on $\overline{\mathbb{IR}}$.

\section{The real vector space $\overline{\mathbb{IR}}$}

\subsection{The semigroup $(\mathbb{IR}$,$+)$}

Consider $x=[x^{-},x^{+}]$ and $y=[y^{-},y^{+}]$ two elements of $\mathbb{IR} $.\ From (\ref{op}) we get the addition
\begin{equation*}
x+y=[x^{-}+y^{-},x^{+}+y^{+}].
\end{equation*}
This operation is commutative, associative and has an unit $[0,0]$ simply denoted by $0$.

\begin{theorem}
The semigroup $(\mathbb{IR},+)$ is commutative and regular.
\end{theorem}

\noindent \textit{Proof. }We recall that a semigroup is a nonempty set with an associative unitary operation $+$. It is
regular if it satisfies
\begin{equation*}
x+z=x+y\Longrightarrow z=y,
\end{equation*}
for all $x,y,z$ . The semigroup $(\mathbb{IR},+\mathbb{)}$ is regular. In fact
\begin{equation*}
x+z=x+y\Longrightarrow \lbrack x^{-}+z^{-},x^{+}+z^{+}]=[x^{-}+y^{-},x^{+}+y^{+}]
\end{equation*}
which gives $z^{-}=y^{-}$ and \ $z^{+}=y^{+},$ that is $z=y.$

\subsection{The group ($\overline{\mathbb{IR}}$,$+)$}

The goal is to define a substraction corresponding to an inverse of the addition. For that we build the symmetrized of
the semigroup ($\mathbb{IR}$,$ +).$ We consider on the set $\mathbb{IR\times IR}$ the equivalence relation:
\begin{equation*}
(x,y)\sim (z,t)\Longleftrightarrow x+t=y+z,
\end{equation*}
for all $x,y,z,t\in \mathbb{IR}$. The quotient set is denoted by $\overline{ \mathbb{IR}}$. The addition of intervals
is compatible with this equivalence relation:
\begin{equation*}
\overline{(x,y)}+\overline{(z,t)}=\overline{(x+z,y+t)}
\end{equation*}
where $\overline{(x,y)}$ is the equivalence class of $(x,y).$ The unit is $ \overline{0}=\{(x,x),x\in \mathbb{IR\}}$
and each element $\overline{(x,y)}$ has an inverse
\begin{equation*}
\smallsetminus \overline{(x,y)}=\overline{(y,x)}.
\end{equation*}

\noindent Then $(\overline{\mathbb{IR}},+)$ is a commutative group.

\medskip

For all $x=[x^{-},x^{+}]\in \mathbb{IR}$, we denote by $l(x)$ his length, so $l(x)=x^{+}-x^{-}$, and by $c(x)$ his
center, so $c(x)=\dfrac{x^{+}+x^{-}}{2} .$

\begin{proposition}
Let $\mathcal{X}=\overline{(x,y)}$ be in $\overline{\mathbb{IR}}$. Thus

\begin{itemize}
\item if $l(y)<l(x),$ there is an unique $A\in \mathbb{IR\setminus R}$ such that $\mathcal{X}=\overline{(A,0)},$

\item if $l(y)>l(x),$ there is an unique $A\in \mathbb{IR}\setminus \mathbb{R }$ such that
    $\mathcal{X}=\overline{(0,A)}=\smallsetminus \overline{(A,0)},$

\item if $l(y)=l(x),$ there is an unique $A=\alpha \in \mathbb{R}$ such that $\mathcal{X}=\overline{(\alpha
    ,0)}=\overline{(0,-\alpha )}. $
\end{itemize}
\end{proposition}

\noindent \textit{Proof. }It is based on the following lemmas:

\begin{lemma}
Consider $\overline{(x,y)}\in \overline{\mathbb{IR}}$ with $l(x)<l(y).$ Then
\begin{equation*}
(x,y)\sim (0,[y^{-}-x^{-},y^{+}-x^{+}]).
\end{equation*}
\end{lemma}

\noindent \textit{Proof.} For all $x,y,z,t\in \mathbb{IR}$, we have
\begin{equation*}
\ (x,y)\sim (z,t)\Longleftrightarrow \left\{
\begin{array}{c}
z^{-}+y^{-}=x^{-}+t^{-}, \\ z^{+}+y^{+}=x^{+}+t^{+}.
\end{array}
\right.
\end{equation*}
If we put $z^{-}=z^{+}=0,$ then
\begin{equation*}
\left\{
\begin{array}{c}
t^{-}=y^{-}-x^{-}, \\ t^{+}=y^{+}-x^{+}.
\end{array}
\right.
\end{equation*}
with the necessary condition $t^{+}>t^{-}$. We obtain $l(x)<l(y).$ So we have
\begin{equation*}
\overline{(x,y)}=\overline{(0,[y^{-}-x^{-},y^{+}-x^{+}])}.
\end{equation*}

\begin{lemma}
Consider $\overline{(x,y)}\in \overline{\mathbb{IR}}$ with $l(y)<l(x),$ then
\begin{equation*}
(x,y)\sim ([x^{-}-y^{-},x^{+}-y^{+}],0).
\end{equation*}
\end{lemma}

\noindent \textit{Proof.} For all $x,y,z\in \mathbb{IR}$, we have
\begin{equation*}
\ (x,y)\sim (z,0)\Longleftrightarrow \left\{
\begin{array}{c}
z^{-}+y^{-}=x^{-} \\ z^{+}+y^{+}=x^{+}
\end{array}
\right. .
\end{equation*}
or
\begin{equation*}
\left\{
\begin{array}{c}
z^{-}=x^{-}-y^{-} \\ z^{+}=x^{+}-y^{+}
\end{array}
\right.
\end{equation*}
with the condition $z^{+}>z^{-}$ which gives $l(y)<l(x).$ So $\overline{(x,y)
}=\overline{([x^{-}-y^{-},x^{+}-y^{+}],0)}.$

\begin{lemma}
Consider $\overline{(x,y)}\in \overline{\mathbb{IR}}$ with $l(x)=l(y),$ then
\begin{equation*}
(x,y)\sim (\alpha ,0)
\end{equation*}
with $\alpha =x^{-}-y^{-}.$
\end{lemma}

\noindent These three lemmas describe the three cases of Proposition 2.

\begin{definition}
Any element $\mathcal{X}=\overline{(A,0)}$ with $A\in \mathbb{IR\setminus R} $ is said positive and we write
$\mathcal{X}>0.$ Any element $\mathcal{X}= \overline{(0,A)}$ with $A\in \mathbb{IR\setminus R}$ is said negative and we
write $\mathcal{X}<0.$ We write $\mathcal{X}\geq \mathcal{X}^{\prime }$ if $ \mathcal{X}\smallsetminus
\mathcal{X}^{\prime }\geq 0.$
\end{definition}

For example if $\mathcal{X}$ \ and $\mathcal{X}^{\prime }$ are positive,
\begin{equation*}
\mathcal{X}\geq \mathcal{X}^{\prime }\Longleftrightarrow l(\mathcal{X})\geq l(\mathcal{X}^{\prime }).
\end{equation*}

\noindent The elements $\overline{(\alpha ,0)}$ with $\alpha \in \mathbb{R} ^{\ast }$ are neither positive nor
negative.

\medskip

\noindent \noindent\textbf{Remark. }This structure of abelian group on $ \overline{\mathbb{IR}}$ has yet been defined
by\ Markov \cite{Ma}. In his paper he presents the Kaucher arithmetic using the completion of the semigroup
$\mathbb{IR}$. Up to now we have sum up this study. In the following section we will develop a topological vectorial
structure. The goal is to define a good differential calculus.

\subsection{ Vector space structure on $\overline{\mathbb{IR}}$}

We are going to construct a real vector space structure on the group $(\overline{\mathbb{IR}},+)$. We recall that if
$A=[a,b]\in \mathbb{IR}$ and $ \alpha \in \mathbb{R}^{+}$, the product $\alpha A$ is the interval $[\alpha a,\alpha
b].$ We consider the external multiplication:
\begin{equation*}
\cdot :\mathbb{R}\times \overline{\mathbb{IR}}\longrightarrow \overline{ \mathbb{IR}}
\end{equation*}
defined, for all $A\in \mathbb{IR}$, by
\begin{equation*}
\left\{
\begin{array}{c}
\alpha \cdot \overline{(A,0)}\text{\ }=\overline{(\alpha A,0)}\text{ ,} \\ \alpha \cdot \overline{(0,A)}\text{\
}=\overline{(0,\alpha A)}\text{ ,}
\end{array}
\right.
\end{equation*}
for all $\alpha >0.$ If $\alpha <0$ we put $\beta =-\alpha $. So we take:
\begin{equation*}
\left\{
\begin{array}{c}
\alpha \cdot \overline{(A,0)}\text{\ }=\overline{(0,\beta A)}, \\ \alpha \cdot \overline{(0,A)}\text{\
}=\overline{(\beta A,0)}.
\end{array}
\right.
\end{equation*}

\noindent We denote $\alpha \mathcal{X}$ instead of $\alpha \cdot \mathcal{X} .$

\begin{lemma}
For any $\alpha \in \mathbb{R}$ and $\mathcal{X}\in \overline{\mathbb{IR}}$ we have:
\begin{equation*}
\left\{
\begin{array}{c}
\alpha (\smallsetminus \mathcal{X})=\smallsetminus (\alpha \mathcal{X}), \\ (-\alpha )\mathcal{X}=\smallsetminus (\alpha
\mathcal{X}).
\end{array}
\right.
\end{equation*}
\end{lemma}

\noindent Indeed, if $\alpha \geq 0$ and $\mathcal{X}=\overline{(A,0)}$ with $A\in \mathbb{IR}$ then $\smallsetminus
\mathcal{X}=\overline{(0,A)}$ and
\begin{equation*}
\alpha (\smallsetminus \mathcal{X})=\overline{(0,\alpha A)}=\smallsetminus \overline{(\alpha A,0)}=\smallsetminus
(\alpha \mathcal{X}).
\end{equation*}
In the same way $(-\alpha )\mathcal{X}=(-\alpha )\overline{(A,0)}=\overline{ (0,\alpha A)}=\smallsetminus (\alpha
\mathcal{X}).$ If $\alpha \leq 0$ and $\mathcal{X}=\overline{(A,0)}$ so $\alpha (\smallsetminus \mathcal{X})=\alpha
\overline{(0,A)}=\overline{(-\alpha A,0)}$ \ and $\smallsetminus (\alpha \mathcal{X})=\smallsetminus
\overline{(0,-\alpha A)}=\overline{(-\alpha A,0)} .$ The calculus is the same for $\mathcal{X}=\overline{(0,A).}$

\begin{proposition}
For all $\alpha ,\beta \in \mathbb{R}$, and for all $\mathcal{X},\mathcal{X} ^{\prime }\in \overline{\mathbb{IR}}$, we
have
\begin{equation*}
\left\{
\begin{tabular}{l}
$(\alpha +\beta )\mathcal{X}=\alpha \mathcal{X}+\beta \mathcal{X}$, \\ $\alpha (\mathcal{X}+\mathcal{X}^{\prime
})=\alpha \mathcal{X}+\alpha \mathcal{X}^{\prime }$, \\ $(\alpha \beta )\mathcal{X}=\alpha (\beta \mathcal{X}).$
\end{tabular}
\right.
\end{equation*}
\end{proposition}

\noindent \textit{Proof}. We are going to study the different cases.

1) The result is trivial when $\alpha ,\beta >0.$

2) Suppose $\alpha >0$ $,$ $\beta <0$ and $\alpha +\beta >0$. We assume $ \mathcal{X}=\overline{(A,0)}$ and $A\in
\mathbb{IR}$. We put $\gamma =-\beta .$ So
\begin{equation*}
(\alpha +\beta )\mathcal{X}=\overline{((\alpha +\beta )A,0)}
\end{equation*}
and
\begin{equation*}
\alpha \mathcal{X}+\beta \mathcal{X}=\overline{(\alpha A,0)}+\overline{ (0,\gamma A)}=\overline{(\alpha A-\gamma
A,0)}=\overline{((\alpha -\gamma )A,0)}=\overline{((\alpha +\beta )A,0)}.
\end{equation*}
The calculus is the same for $\mathcal{X}=\overline{(0,A)}.$ So we have $ (\alpha +\beta )\mathcal{X}=\alpha
\mathcal{X}+\beta \mathcal{X}$ for all $\mathcal{X}\in \overline{\mathbb{IR}}.$

3) If $\alpha >0$ $,$ $\beta <0$ $\ $\ and $\alpha +\beta <0$, the proof is the same that the previous case.

4) If $\beta >0$ and $\alpha <0,$ we refer to the cases $2)$ et $3)$.

5) If $\alpha ,\beta <0$, we put $\beta _{1}=-\beta $ and $\alpha _{1}=-\alpha $. Thus $\alpha _{1},\beta _{1}>0,$ and
we find again the first case.

\bigskip

\noindent So we have the result:

\begin{theorem}
The triplet $(\overline{\mathbb{IR}},+,\cdot )$ is a real vector space.
\end{theorem}

\subsection{Basis and dimension}

We consider the vectors $\mathcal{X}_{1}=\overline{([0,1],0)}$ and $\mathcal{ X}_{2}=\overline{([1,1],0)}$ of
$\overline{\mathbb{IR}}.$

\begin{theorem}
The family $\{\mathcal{X}_{1},\mathcal{X}_{2}\}$ is a basis of $\overline{\mathbb{IR}}.$ So $\dim
_{\mathbb{R}}\overline{\mathbb{IR}}=2.$
\end{theorem}

\noindent \textit{Proof. }We have the following decompositions:\begin{equation*}
\left\{
\begin{array}{c}
\overline{([a,b],0)\text{ }}=(b-a)\mathcal{X}_{1}+a\mathcal{X}_{2}, \\ \overline{(0,[c,d])\text{
}}=(c-d)\mathcal{X}_{1}-c\mathcal{X}_{2}.
\end{array}
\right.
\end{equation*}
The linear map
\begin{equation*}
\varphi :\overline{\mathbb{IR}}\longrightarrow \mathbb{R}^{2}
\end{equation*}
defined by
\begin{equation*}
\left\{
\begin{array}{l}
\varphi (\,\overline{([a,b],0)}\text{ })=(b-a,a), \\ \varphi (\,\overline{(0,[c,d])}\text{ })=(c-d,-c)
\end{array}
\right.
\end{equation*}
is a linear isomorphism and $\overline{\mathbb{IR}}$ is canonically isomorphic to $\mathbb{R}^{2}$.

\medskip

\noindent \textbf{Remark. }Let $E$ be the subspace generated by $\mathcal{X} _{2}.$ The vectors of $E$ correspond to
the elements which have a non defined sign. Then the relation $\leq $ defined in the paragraph $1.2$ gives an order
relation on the quotient space $\overline{\mathbb{IR}}/E.$

\subsection{A Banach structure on $\overline{\mathbb{IR}}$}

Let us begin to define a norm on $\overline{\mathbb{IR}}.$ Any element $ \mathcal{X}\in \overline{\mathbb{IR}}$ $\ $is
written $\overline{(A,0)}$ or $ \overline{(0,A)}.$\ We define its length $l(\mathcal{X})$ as the length of $ A $ and
its center as $c(A)$ or $-c(A)$ in the second case.

\begin{theorem}
The map $||$ $||$ $:$ $\overline{\mathbb{IR}}\longrightarrow $ $\mathbb{R}$ given by
\begin{equation*}
||\mathcal{X}||=l(\mathcal{X})+|c(\mathcal{X})|
\end{equation*}
for any $\mathcal{X}\in $ $\overline{\mathbb{IR}}$ is a norm$.$
\end{theorem}

\noindent \noindent \noindent \textit{Proof}. We have to verify the following axioms:
\begin{equation*}
\left\{
\begin{array}{l}
1)\text{ }||\mathcal{X}||=0\Longleftrightarrow \mathcal{X}=0, \\ 2)\text{ }\forall \lambda \in \mathbb{R}\text{
}||\lambda \mathcal{X} ||=|\lambda |||\mathcal{X}||,\text{ } \\ 3)\text{ }||\mathcal{X}+\mathcal{X}^{\prime }||\leq
||\mathcal{X}||+|| \mathcal{X}^{\prime }||.
\end{array}
\right.
\end{equation*}

\noindent 1) If $||\mathcal{X}||=0$, then $l(\mathcal{X})=|c(\mathcal{X})|=0$ and $\mathcal{X}=0.$

\bigskip

\noindent \noindent \noindent \noindent 2) Let $\lambda \in \mathbb{R}.$ We have
\begin{equation*}
||\lambda \mathcal{X}||=l(\lambda \mathcal{X})+|c(\lambda \mathcal{X} )|=|\lambda |l(\mathcal{X})+|\lambda
||c(\mathcal{X})|=|\lambda |||\mathcal{X }||.
\end{equation*}

\medskip

\noindent \noindent 3) We consider that $I$ refers to $\mathcal{X}$ and $J$ refers to $\mathcal{X}^{\prime}$ thus
$\mathcal{X}=\overline{(I,0)}$ or $= \overline{(0,I)}$. We have to study the two different cases:

\noindent i) If $\mathcal{X}+\mathcal{X}^{\prime }=\overline{(I+J,0)}$ or $ \overline{(0,I+J)}$, then
\begin{eqnarray*}
||\mathcal{X}+\mathcal{X}^{\prime }|| &=&l(I+J)+|c(I+J)|=l(I)+l(J)+|c(I)+c(J)|\leq l(I)+|c(I)|+l(J)+|c(J)| \\
&=&||\mathcal{X}||+||\mathcal{X}^{\prime }||.
\end{eqnarray*}

\noindent ii) Let $\mathcal{X}+\mathcal{X}^{\prime }=\overline{(I,J)}.$ If $ \overline{(I,J)}=\overline{(K,0)}$ then
$K+J=I$ and
\begin{equation*}
||\mathcal{X}+\mathcal{X}^{\prime }||=||\overline{(K,0)} ||=l(K)+|c(K)|=l(I)-l(J)+|c(I)-c(J)|
\end{equation*}
that is
\begin{equation*}
||\mathcal{X}+\mathcal{X}^{\prime }||\leq l(I)+|c(I)|-l(J)+|c(J)|\leq
l(I)+|c(I)|+l(J)+|c(J)|=||\mathcal{X}||+||\mathcal{X}^{\prime }||.
\end{equation*}
So we have a norm on $\overline{\mathbb{IR}}.$

Now we shall show that $\overline{\mathbb{IR}}$ is a Banach space that is all Cauchy sequences converge in
$\overline{\mathbb{IR}}$.

\begin{theorem}
The normed vector space $\overline{\mathbb{IR}}$ is a Banach space.
\end{theorem}

\noindent \textit{Proof}. We recall that a sequence $(\mathcal{X}_{n})_{n\in \mathbb{N}}$ with $\mathcal{X}_{n}\in
\overline{\mathbb{IR}}$ is a Cauchy sequence if
\begin{equation*}
\forall \varepsilon >0,\exists N\in \mathbb{N},\forall n,m\geq N,||\mathcal{X }_{n}\smallsetminus \mathcal{X}_{m}||\leq
\varepsilon .
\end{equation*}
We have to verify that all Cauchy sequences in $\overline{\mathbb{IR}}$ converge in $\overline{\mathbb{IR}}$. Let
$(\mathcal{X}_{n})_{n\in \mathbb{N} }$ be a Cauchy sequence. We suppose that for all $n,m\geq N$ \ we have $
\mathcal{X}_{n}=\overline{(A_{n},0)}$ and $\mathcal{X}_{m}=\overline{ (A_{m},0)}$ . So we have
$\mathcal{X}_{n}\smallsetminus \mathcal{X}_{m}= \overline{(A_{n,}A_{m})}$ and
\begin{equation*}
||\mathcal{X}_{n}\smallsetminus \mathcal{X} _{m}||=|l(A_{n})-l(A_{m})|+|c(A_{n})-c(A_{m})|\leq \varepsilon .
\end{equation*}
Then for $n,m>N,$ $|l(A_{n})-l(A_{m})|\leq \varepsilon $ and $ |c(A_{n})-c(A_{m})|\leq \varepsilon .$ The sequences
$(l(A_{n}))_{n\in \mathbb{N}}$\ and $(c(A_{n}))_{n\in \mathbb{N}}$ are Cauchy sequences in $ \mathbb{R}$, which is a
complete space, thus the two sequences converge. Let $l$ and $c$ be their limits. It exists an unique $L$ so that
$l(L)=l$ and $ c(L)=c$. Thus we can thus deduce that $(\mathcal{X}_{n})_{n\in \mathbb{N}}$ converges to $L$. It is the
same if all terms of the sequence are of type $ \mathcal{X}_{n}=\overline{(0,A_{n})}.$ We suppose now that the sequence
$( \mathcal{X}_{n})_{n\in \mathbb{N}}$ is not of constant sign from a certain rank. So
\begin{equation*}
\exists n,m>N\text{ such that }\mathcal{X}_{n}=\overline{(A_{n},0)}\text{ and }\mathcal{X}_{m}=\overline{(0,A_{m})},
\end{equation*}
and
\begin{equation*}
||\mathcal{X}_{n}\smallsetminus \mathcal{X}_{m}||=||\overline{(A_{n}+A_{m},0)
}||=l(A_{n})+l(A_{m})+|c(A_{n})+c(A_{m})|\leq \varepsilon .\qquad (\ast )
\end{equation*}
We deduce that $l(A_{n})\underset{n\rightarrow +\infty }{\longrightarrow }0.$ Moreover, we can consider the subsequence
whose terms are of type $( \overline{\text{ }(A_{n},0)})_{n\in \mathbb{N}}$. It converges to $\mathcal{X }$ with
$l(\mathcal{X})=0$ thus $\mathcal{X}=\overline{(\alpha ,0)}$. In the same way we can consider the subsequence of terms
are $(\overline{\text{ } (0,A_{n})})$. It converges to $\mathcal{X}^{\prime }=\overline{(0,\beta )}$ since
$l(\mathcal{X}^{\prime })=0$. But the Inequation $(\ast )$ implies $ \alpha +\beta \leq \varepsilon $ for all
$\varepsilon $ thus $\mathcal{X} ^{\prime }=$ $\overline{(0,-\alpha )}=\mathcal{X}.$ It follows that $
\overline{\mathbb{IR}}$ is complete.

\subsection{Description of a $\protect\varepsilon -$neighbourhood of $
\mathcal{X}_{0}\in $ $\overline{\mathbb{IR}}$}

We assume that $\mathcal{X}_{0}=\overline{([a,b],0)}$ with $0<a<b.$ Let $\mathcal{X}=\overline{([x,y],0)}.$ Then
$\mathcal{ X\smallsetminus X}_{0}=\overline{([x,y],[a,b])}$.

\begin{enumerate}
\item $x-a<y-b$ that is $l(\mathcal{X}_{0})<l(\mathcal{X}).$ Then $\mathcal{ X\smallsetminus
    X}_{0}=\overline{([x-a,y-b],0)}$ and
$$||\mathcal{
X\smallsetminus X}_{0}||=y-b-x+a+|\dfrac{y-b+x-a}{2}|.$$

\begin{itemize}
\item $y+x>b+a.$ We have $||\mathcal{X\smallsetminus X}_{0}||=\dfrac{3y-x-3b+a}{2}$ and
    $||\mathcal{X\smallsetminus X}_{0}||<\varepsilon $ is represented by the convex set defined by
\begin{equation*}
\left\{
\begin{array}{c}
y>x+b-a, \\ y>-x+b+a,
\end{array}
\right.
\end{equation*}
that
\begin{equation*}
x+3b-a-2\varepsilon <3y<x+3b-a+2\varepsilon .
\end{equation*}
It is a triangle whose vertices are $A(a-\dfrac{\varepsilon }{2},b+\dfrac{ \varepsilon }{2}),$ $B(a+\varepsilon
,b+\varepsilon ),$ $C(a,b).$

\item $y+x>b+a.$ We have $||\mathcal{X\smallsetminus X}_{0}||=\dfrac{ y-3x-b+3a}{2}$ and
    $||\mathcal{X\smallsetminus X}_{0}||<\varepsilon $ is represented by the convex set defined by
\begin{equation*}
\left\{
\begin{array}{c}
y>x+b-a, \\ y<-x+b+a,
\end{array}
\right.
\end{equation*}
that
\begin{equation*}
3x+3b-3a-2\varepsilon <y<3x+b-3a+2\varepsilon .
\end{equation*}
It is a triangle whose vertices are $A,$ $B^{\prime }(a-\varepsilon ,b-\varepsilon ),$ $C.$
\end{itemize}

\item $x-a>y-b$ that is $l(\mathcal{X}_{0})<l(\mathcal{X}).$ Then $\mathcal{ X\smallsetminus
    X}_{0}=\overline{(0,[a-x,b-y])}$ and $$||\mathcal{ X\smallsetminus X}_{0}||=b-y-a+x+|\dfrac{b-y+a-x}{2}|.$$

\begin{itemize}
\item $y+x>b+a,$ we have $||\mathcal{X\smallsetminus X}_{0}||=\dfrac{ 3x-y+b-3a}{2}$ and
    $||\mathcal{X\smallsetminus X}_{0}||<\varepsilon $ is represented by the convex set defined by
\begin{equation*}
\left\{
\begin{array}{c}
y<x+b-a, \\ y>-x+b+a,
\end{array}
\right.
\end{equation*}
that
\begin{equation*}
-3x-b+3a-2\varepsilon <-y<-3x-b+3a+2\varepsilon .
\end{equation*}
It is a triangle whose vertices are $B,$ $B^{\prime \prime }(a+\dfrac{ \varepsilon }{2},b-\dfrac{\varepsilon
}{2}),$ $C.$

\item $y+x>b+a,$ we have $||\mathcal{X\smallsetminus X}_{0}||=\dfrac{ x-3y+3b-a}{2}$ and
    $||\mathcal{X\smallsetminus X}_{0}||<\varepsilon $ is represented by the triangle whose vertices are
    $B^{\prime },$ $B^{\prime \prime },$ $C.$
\end{itemize}
\end{enumerate}

\noindent If $\mathcal{X}=\overline{(0,[x,y])},$ then $\mathcal{X\smallsetminus X}_{0}=\overline{(0,[x+a,y+b])}$ and
$$||\mathcal{X\smallsetminus X}_{0}||=
\dfrac{3y-x+3b-a}{2}=||\mathcal{X}||\mathcal{+}||\mathcal{X}_{0}||\nless \varepsilon. $$ Such a point cannot be in a
neighborhood of $\mathcal{X}_{0}$.

\begin{proposition}
An $\varepsilon$-neighborhood of $\mathcal{X}_{0}=\overline{([a,b],0)}$ with $0\leq a\leq b$ is represented by the
parallelogram  whose vertices are
$A_1=(a-\varepsilon,b-\varepsilon),A_2=(a+\frac{\varepsilon}{2},b-\frac{\varepsilon}{2}),A_3=(a+\varepsilon,b+\varepsilon),A_4=(a-\frac{\varepsilon}{2},b+\frac{\varepsilon}{2})$.
\end{proposition}

\section{\protect\bigskip A 4-dimensional associative algebra associated to $
\overline{\mathbb{IR}}$}

\subsection{Classical product of intervals}

We consider $X$,$Y\in \mathbb{IR}$. The multiplication of intervals is defined by
\begin{equation*}
X\cdot Y=[\min (x^{-}y^{-},x^{-}y^{+},x^{+}y^{-},x^{+}y^{+}),\max (x^{-}y^{-},x^{-}y^{+},x^{+}y^{-},x^{+}y^{+})].
\end{equation*}
Let $\mathcal{X}=\overline{(X,0)}$ and $\mathcal{X}^{\prime }=\overline{(Y,0) \text{ }}$ be in
$\overline{\mathbb{IR}}.$ We put
\begin{equation*}
\mathcal{XX}^{\prime }=\overline{(XY,0)}.
\end{equation*}
For this product we have:

\begin{proposition}
For all $\mathcal{X}=\overline{(X,0)}$ and $\mathcal{X}^{\prime }=\overline{ (Y,0)}$ in $\overline{\mathbb{IR}},$ we
have
\begin{equation*}
||\mathcal{XX}^{\prime }||\leq ||\mathcal{X}||\text{ }||\mathcal{X}^{\prime }||.
\end{equation*}
\end{proposition}

\noindent\textit{Proof. }In the following table, the boxes represent $||\mathcal{XX}^{\prime }||$ following the values
of $||\mathcal{X}||$ and $|| \mathcal{X}^{\prime }||.$

\noindent If $\mathcal{X}=\overline{([x_{1},x_{2}],0)}$ then
\begin{equation*}
\left\{
\begin{array}{l}
||\mathcal{X}||=\dfrac{3x_{2}-x_{1}}{2}\text{ if \ }c(\mathcal{X})>0, \\ ||\mathcal{X}||=\dfrac{x_{2}-3x_{1}}{2}\text{
if }c(\mathcal{X})<0.
\end{array}
\right.
\end{equation*}

\noindent Considering the different situations, we obtain
\begin{equation*}
||\mathcal{X}||||\mathcal{X}^{\prime }||-||\mathcal{XX}^{\prime }||=\frac{3}{ 4}l(\mathcal{X})l(\mathcal{X}^{\prime })
\end{equation*}
or $\frac{1}{2}||\mathcal{X}||l(\mathcal{X}^{\prime })$ or $\frac{1}{2}|| \mathcal{X}^{\prime }||l(\mathcal{X}).$ These
expressions are always positive. We have $||\mathcal{X}||||\mathcal{X}^{\prime }||=||\mathcal{XX} ^{\prime }||$ if
$\mathcal{X}$ or $\mathcal{X}^{\prime }$ are reduce to one point.

\begin{proposition}
We consider $\mathcal{X}=\overline{(X,0)}$ and $\mathcal{X}^{\prime }= \overline{(Y,0)}$ in $\overline{\mathbb{IR}}$.
We have
\begin{equation*}
X\subset Y\Rightarrow ||\mathcal{X}||\leq ||\mathcal{X}^{\prime }||.
\end{equation*}
\end{proposition}

\noindent \noindent \noindent \textit{Proof. }Consider $X=[x_{1},x_{2}]$ and $Y=[y_{1},y_{2}].$

\noindent \textit{First case: }$y_{1}\geq 0$. So $2||\mathcal{X}^{\prime }||=3y_{2}-y_{1}.$ As $X\subset Y$, then
$2||\mathcal{X}||=3x_{2}-x_{1}$ and $||\mathcal{X}||\leq ||\mathcal{X}^{\prime }||.$

\noindent \textit{Second case: }$y_{1}<0,y_{2}>0$. If $c(Y)\geq 0,$ so $2|| \mathcal{X}^{\prime }||=3y_{2}-y_{1}.$ If
$c(X)\geq 0$, from the first case $ ||\mathcal{X}||\leq ||\mathcal{X}^{\prime }||.$ Otherwise $2||\mathcal{X}
||=x_{2}-3x_{1}.$ Thus $||\mathcal{X}||\leq ||\mathcal{X}^{\prime }||$ if and only if $3y_{2}-y_{1}\geq x_{2}-3x_{1},$
that is $3(y_{2}+x_{1})\geq x_{2}+y_{1}$ which is true.

\noindent If $c(Y)\leq 0$, then $2||\mathcal{X}^{\prime }||=y_{2}-3y_{1}.$ If $c(X)\leq 0$, thus
$2||\mathcal{X}||=x_{2}-3x_{1}$ and $||\mathcal{X} ||\leq ||\mathcal{X}^{\prime }||$. If $c(X)\geq 0$,
$||\mathcal{X}||\leq || \mathcal{X}^{\prime }||$ is equivalent to $y_{2}-3y_{1}\geq 3x_{2}-x_{1}.$ But $c(Y)\leq 0$
implies $y_{1}+y_{2}\leq 0$ and $y_{2}-3y_{1}\geq 4y_{2}$. Similarly $3x_{2}-x_{1}\leq 4x_{2},$ thus $y_{2}-3y_{1}\geq
3x_{2}-x_{1}$ because $x_{2}\leq y_{2}.$

\noindent\textit{Third case: }$y_{1}<0,$ $y_{2}<0.$ Similar computations give the result.\bigskip

\noindent \noindent \noindent \textbf{Remark. }If $\mathcal{X}>0,$ i.e $ \mathcal{X}=\overline{(X,0)},$ and
$\mathcal{X}^{\prime }<0,$ i.e. $\mathcal{ X}^{\prime }=\overline{(0,Y)},$ so $\smallsetminus \mathcal{X}^{\prime }>0$
and if $X\subset Y$ we deduce $||\mathcal{X}||\leq ||\smallsetminus \mathcal{ X}^{\prime }||=||\mathcal{X}^{\prime
}||.$

\bigskip

\subsection{Definition of $\mathcal{A}_{4}$}

Recall that by an algebra we mean a real vector space with an associative ring structure. Consider the $4$-dimensional
associative algebra whose product in a basis $\{e_{1},e_{2},e_{3},e_{4}\}$ is given by
\begin{equation*}
\begin{tabular}{|l|l|l|l|l|}
\hline & $e_{1}$ & $e_{2}$ & $e_{3}$ & $e_{4}$ \\ \hline $e_{1}$ & $e_{1}$ & $0$ & $0$ & $e_{4}$ \\ \hline $e_{2}$ & $0$
& $e_{2}$ & $e_{3}$ & $0$ \\ \hline $e_{3}$ & $0$ & $e_{3}$ & $e_{2}$ & $0$ \\ \hline $e_{4}$ & $e_{4}$ & $0$ & $0$ &
$e_{1}$ \\ \hline
\end{tabular}
.
\end{equation*}
The unit is the vector $e_{1}+e_{2}.$ This algebra is a direct sum of two ideals: $\mathcal{A}_{4}=I_{1}+I_{2}$ where
$I_{1}$ is generated by $e_{1}$ and $e_{4}$ and $I_{2}$ is generated by $e_{2}$ and $e_{3}.$ It is not an integral
domain, that is, we have divisors of $0.$ For example $e_{1}\cdot e_{2}=0.$

\begin{proposition}
\bigskip The ring $\mathcal{A}_{4}$ is principal that is every ideal is
generated by one element.
\end{proposition}

\noindent \textit{Proof. }In fact we have only two ideals $I_{1}$ and $I_{2}$ and $I_{1}$ is generated by $e_{4}$ and
$I_{2}$ is generated by $e_{3}.$ We denote by $\mathcal{A}_{4}^{\ast }$ the group of invertible elements. We compute
this group. The cartesian expression of this product is, for $ x=(x_{1},x_{2},x_{3},x_{4})$ and
$y=(y_{1},y_{2},y_{3},y_{4})$ in $\mathcal{A }_{4}$:
\begin{equation*}
x\cdot y=(x_{1}y_{1}+x_{4}y_{4},x_{2}y_{2}+x_{3}y_{3},x_{3}y_{2}+x_{2}y_{3},x_{4}y_{1}+x_{1}y_{4}).
\end{equation*}
We consider the equation
\begin{equation*}
x\cdot y=(1,1,0,0).
\end{equation*}
We obtain
\begin{equation*}
\left\{
\begin{array}{c}
x_{1}y_{1}+x_{4}y_{4}=1, \\ x_{2}y_{2}+x_{3}y_{3}=1, \\ x_{3}y_{2}+x_{2}y_{3}=0, \\ x_{4}y_{1}+x_{1}y_{4}=0.
\end{array}
\right.
\end{equation*}
For a given vector $x$, we obtain a solution $y$ if and only if:
\begin{equation*}
(x_{1}^{2}-x_{4}^{2})(x_{2}^{2}-x_{3}^{2})\neq 0.
\end{equation*}

\begin{proposition}
The multiplicative group $\mathcal{A}_{4}^{\ast }$ \ is the set of elements $ x=(x_{1},x_{2},x_{3},x_{4})$ such that
\begin{equation*}
\left\{
\begin{array}{c}
x_{4}\neq \pm x_{1}, \\ x_{3}\neq \pm x_{2}.
\end{array}
\right.
\end{equation*}
If $x\in $ $\mathcal{A}_{4}^{\ast }$ we have:
\begin{equation*}
x^{-1}=\left( \frac{x_{1}}{x_{1}^{2}-x_{4}^{2}},\frac{x_{2}}{
x_{2}^{2}-x_{3}^{2}},\frac{x_{3}}{x_{2}^{2}-x_{3}^{2}},\frac{x_{4}}{ x_{1}^{2}-x_{4}^{2}}\right) .
\end{equation*}
\end{proposition}

\subsection{An embedding of $\mathbb{IR}$ in $\mathcal{A}_{4}$}

We define a correspondence between $\mathbb{IR}$ and $\mathcal{A}_{4}.$ Let $ \varphi $ be the map
\begin{equation*}
\varphi :\mathbb{IR\longrightarrow }\mathcal{A}_{4}
\end{equation*}
defined by:
\begin{equation*}
\varphi (X)=\left\{
\begin{array}{l}
(x_{1},x_{2},0,0)\text{ if }x_{1},x_{2}\geq 0, \\ (0,x_{2},-x_{1},0)\text{ if }x_{1}\leq 0\text{ and }x_{2}\geq 0, \\
(0,0,-x_{1},-x_{2})\text{ if }x_{1},x_{2}\leq 0,
\end{array}
\right.
\end{equation*}

for every $X=[x_{1},x_{2}]$ in $\mathbb{IR}$.

\begin{theorem}
Let $X=[x_{1},x_{2}]$ and $Y=[y_{1},y_{2}]$ be in $\mathbb{IR}$.

\begin{itemize}
\item If $x_{1}x_{2}>0$ or $y_{1}y_{2}>0,$ then
\begin{equation*}
\varphi (XY)=\varphi (X)\cdot \varphi (Y).
\end{equation*}

\item If $x_{1}x_{2}<0$ and $y_{1}y_{2}<0$, then
\begin{equation*}
\varphi ^{-1}(\varphi (XY))\subset \varphi ^{-1}(\varphi (X)\cdot \varphi (Y)).
\end{equation*}
\end{itemize}
\end{theorem}

\noindent \textit{Proof.} For the first case we have the following table:
\begin{equation*}
\begin{tabular}{|l|l|l|l|}
\hline & $\varphi (X)$ & $\varphi (Y)$ & $\varphi (XY)$ \\ \hline $\left\{
\begin{array}{c}
x_{1},x_{2}>0 \\ y_{1},y_{2}>0
\end{array}
\right. $ & $(x_{1},x_{2},0,0)$ & $(y_{1},y_{2},0,0)$ & $ (x_{1}y_{1},x_{2}y_{2},0,0)$ \\ \hline $\left\{
\begin{array}{l}
x_{1},x_{2}>0 \\ y_{1}<0,y_{2}>0
\end{array}
\right. $ & $(x_{1},x_{2},0,0)$ & $(0,y_{2},-y_{1},0)$ & $ (0,x_{2}y_{2},-x_{2}y_{1},0)$ \\ \hline $\left\{
\begin{array}{c}
x_{1},x_{2}>0 \\ y_{1},y_{2}<0
\end{array}
\right. $ & $(x_{1},x_{2},0,0)$ & $(0,0,-y_{1}-,y_{2})$ & $ (0,0,-x_{2}y_{1},-x_{1}y_{2})$ \\ \hline $\left\{
\begin{array}{l}
x_{1},x_{2}<0 \\ y_{1}<0,y_{2}>0
\end{array}
\right. $ & $(0,0,-x_{1},-x_{2})$ & $(0,y_{2},-y_{1},0)$ & $ (0,x_{1}y_{1},-x_{1}y_{2},0)$ \\ \hline $\left\{
\begin{array}{c}
x_{1},x_{2}<0 \\ y_{1,}y_{2}<0
\end{array}
\right. $ & $(0,0,-x_{1},-x_{2})$ & $(0,0,-y_{1}-,y_{2})$ & $ (x_{2}y_{2},x_{1}y_{1},0,0)$ \\ \hline
\end{tabular}
\end{equation*}
and we see that in each case we have $\varphi (XY)=\varphi (X)\cdot \varphi (Y).$ In the second case,
\begin{equation*}
\varphi (X)=(0,x_{2},-x_{1},0)\text{ and }\varphi (Y)=(0,y_{2},-y_{1},0) \text{ . }
\end{equation*}
Then $\varphi (X)\cdot \varphi (Y)=(0,x_{2}y_{2}+x_{1}y_{1},-x_{1}y_{2}-x_{2}y_{1},0).$ But $XY$ can be equal to one of
the following intervals
\begin{equation*}
\begin{tabular}{ll}
$XY=$ & $\left\{
\begin{array}{c}
\lbrack x_{1}y_{2},x_{2}y_{2}], \\ \lbrack x_{2}y_{1},x_{2}y_{2}], \\ \lbrack x_{1}y_{2},x_{1}y_{1}], \\ \lbrack
x_{2}y_{1},x_{1}y_{1}].
\end{array}
\right. $
\end{tabular}
\end{equation*}
Then
\begin{equation*}
\varphi (XY)\in \left\{
(0,x_{2}y_{2},-x_{1}y_{2},0),(0,x_{2}y_{2},-x_{2}y_{1},0),(0,x_{1}y_{1},-x_{1}y_{2},0),(0,x_{1}y_{1},-x_{2}y_{1},0)\right\}
\end{equation*}
and
\begin{equation*}
\varphi ^{-1}(\varphi (X)\cdot \varphi (Y))=[x_{1}y_{2}+x_{2}y_{1},x_{2}y_{2}+x_{1}y_{1}]\supset \varphi ^{-1}(\varphi
(XY)).
\end{equation*}

\bigskip

\noindent\textbf{Remark. }This inclusion is, in some sense, universal because it contains all the cases and it is the
minimal expression satisfying this property. For a computational use of this product, this is not a problem because the
result $\varphi (X)\cdot \varphi (Y)$ contains the classical product (it is a constraint asked by the arithmetic
programming).

\subsection{An embedding of $\overline{\mathbb{IR}}$ in $\mathcal{A}_{4}$}

>From the map
\begin{equation*}
\varphi :\mathbb{IR\longrightarrow }\mathcal{A}_{4}
\end{equation*}
we would like to define $\overline{\varphi }:$ $\overline{\mathbb{IR}} \mathbb{\longrightarrow }\mathcal{A}_{4}$ such
that
\begin{equation*}
\overline{\varphi }\overline{(K,0)}=\varphi (K).
\end{equation*}
To define $\overline{\varphi }\overline{(0,K)}$, we put
\begin{equation*}
\overline{\varphi }\overline{(0,K)}=-\overline{\varphi }\overline{(K,0)}.
\end{equation*}
If we consider $K=[x_{1},x_{2}]$, then we have
\begin{equation*}
\begin{tabular}{lll}
$\overline{\varphi }\overline{(0,K)}$ & $=$ & $\left\{
\begin{array}{l}
-x_{1},-x_{2},0,0)\text{ if }x_{1}\geq 0, \\ (0,-x_{2},-x_{1},0)\text{ if }x_{1}x_{2}\leq 0, \\
(0,0,-x_{1},-x_{2})\text{ if }x_{2}\leq 0.
\end{array}
\right. $
\end{tabular}
\end{equation*}
Thus the image of $\overline{\mathbb{IR}}$ in $\mathcal{A}_{4}$ is constituted of the elements
\begin{equation*}
\left\{
\begin{array}{l}
(x_{1},x_{2},0,0)\text{ with }0\leq x_{1}\leq x_{2}\text{ which corresponds to }([x_{1},x_{2}],0), \\
(0,x_{2},-x_{1},0)\text{ with }x_{1}\leq 0\leq x_{2}\text{ which corresponds to }([x_{1},x_{2}],0), \\
(0,0,-x_{1},-x_{2})\text{ with }x_{1}\leq x_{2}\leq 0\text{ which corresponds to }([x_{1},x_{2}],0), \\
(-x_{1},-x_{2},0,0)\text{with }0\leq x_{1}\leq x_{2}\text{ which corresponds to }(0,[x_{1},x_{2}]), \\
(0,-x_{2},x_{1},0)\text{ with }x_{1}\leq 0\leq x_{2}\text{ which corresponds to }(0,[x_{1},x_{2}]), \\
(0,0,x_{1},x_{2})\text{ with }x_{1}\leq x_{2}\leq 0\text{ which corresponds to }(0,[x_{1},x_{2}]).
\end{array}
\right.
\end{equation*}
The map $\overline{\varphi }:$ $\overline{\mathbb{IR}}\mathbb{ \longrightarrow }\mathcal{A}_{4}$ \ is not linear. For
example, if $\mathcal{ X}_{1}=\overline{([2,4],0)}$ and $\mathcal{X}_{2}=\overline{(0,[1,6])}$, then
\begin{equation*}
\mathcal{X}_{1}+\mathcal{X}_{2}=\overline{([2,4],[1,6])}=\overline{(0,[-1,2]) },
\end{equation*}
and
\begin{equation*}
\overline{\varphi }(\mathcal{X}_{1}+\mathcal{X}_{2})=(0,-2,-1,0).
\end{equation*}
Moreover
\begin{equation*}
\overline{\varphi }(\mathcal{X}_{1})+\overline{\varphi }(\mathcal{X} _{2})=(2,4,0,0)+(-1,-6,0,0)=(1,-2,0,0),
\end{equation*}
and $(1,-2,0,0)$ is not in the image of $\overline{\varphi }$. Let us introduce in $\mathcal{A}_{4}$ the following
equivalence relation $\mathcal{R }$ given by
\begin{equation*}
(x_{1},x_{2},x_{3},x_{4})\sim (y_{1},y_{2},y_{3},y_{4})\Longleftrightarrow \left\{
\begin{array}{c}
x_{1}-y_{1}=x_{3}-y_{3}, \\ x_{2}-y_{2}=x_{4}-y_{4}
\end{array}
\right.
\end{equation*}
and consider the map
\begin{equation*}
\overline{\overline{\varphi }}:\overline{\mathbb{IR}}\longrightarrow
\overline{\mathcal{A}_{4}}=\frac{\mathcal{A}_{4}}{\mathcal{R}}
\end{equation*}
given by $\overline{\overline{\varphi }}=\Pi \circ \overline{\varphi }$ where $\Pi $ is a canonical projection. This map
is surjective. In fact we have the correspondence

\begin{itemize}
\item $x_{1}-x_{3}\geq 0$, $x_{2}-x_{4}\geq 0$, $x_{1}-x_{3}\leq x_{2}-x_{4}$
\begin{equation*}
(x_{1},x_{2},x_{3},x_{4})\sim (x_{1}-x_{3},x_{2}-x_{4},0,0)=\overline{ \varphi }([x_{1}-x_{3},x_{2}-x_{4}],0).
\end{equation*}

\item $x_{1}-x_{3}\geq 0$, $x_{2}-x_{4}\geq 0$, $x_{1}-x_{3}\geq x_{2}-x_{4}$
\begin{equation*}
(x_{1},x_{2},x_{3},x_{4})\sim (0,0,x_{3}-x_{1},x_{4}-x_{2})=\overline{ \varphi }(0,[x_{3}-x_{1},x_{4}-x_{2}]).
\end{equation*}

\item $x_{1}-x_{3}\geq 0$, $x_{2}-x_{4}\leq 0$
\begin{equation*}
(x_{1},x_{2},x_{3},x_{4})\sim (0,x_{2}-x_{4},x_{3}-x_{1},,0)=\overline{ \varphi }(0,[x_{3}-x_{1},x_{4}-x_{2}]).
\end{equation*}

\item $x_{1}-x_{3}\leq 0$, $x_{2}-x_{4}\geq 0$
\begin{equation*}
(x_{1},x_{2},x_{3},x_{4})\sim (0,x_{2}-x_{4},x_{3}-x_{1},,0)=\overline{ \varphi }([x_{3}-x_{1},x_{2}-x_{4}],0).
\end{equation*}

\item $x_{1}-x_{3}\leq 0$, $x_{2}-x_{4}\leq 0$, $x_{1}-x_{3}\geq x_{2}-x_{4}$
\begin{equation*}
(x_{1},x_{2},x_{3},x_{4})\sim (x_{1}-x_{3},x_{2}-x_{4},0,0)=\overline{ \varphi }(0,[x_{3}-x_{1},x_{4}-x_{2}]).
\end{equation*}

\item $x_{1}-x_{3}\leq 0$, $x_{2}-x_{4}\leq 0$, $x_{1}-x_{3}\leq x_{2}-x_{4}(x_{1},x_{2},x_{3},x_{4})\sim
    (0,0,x_{3}-x_{1},x_{4}-x_{2})= \overline{\varphi }([x_{3}-x_{1},x_{2}-x_{4}],0).$
\end{itemize}

\noindent This correspondence defines a map
\begin{equation*}
\psi :\overline{\mathcal{A}_{4}}\longrightarrow \overline{\mathbb{IR}}.
\end{equation*}
In the following, to simplifies notation, we write $\overline{\varphi }$ instead of $\overline{\overline{\varphi }}.$

\begin{definition}
For any $\mathcal{X},\mathcal{X}^{\prime }\in \overline{\mathbb{IR}}$, we put
\begin{equation*}
\mathcal{X}\bullet \mathcal{X}^{\prime }=\psi (\overline{\overline{\varphi }( \mathcal{X})}\bullet
\overline{\overline{\varphi }(\mathcal{X}^{\prime })}).
\end{equation*}
\end{definition}

This multiplication is distributive with respect the the addition. In fact
\begin{equation*}
(\mathcal{X}_{1}+\mathcal{X}_{2})\bullet \mathcal{X}^{\prime }=\psi ( \overline{\overline{\varphi
}(\mathcal{X}_{1}+\mathcal{X}_{2})}\bullet \overline{\overline{\varphi }(\mathcal{X}^{\prime })}).
\end{equation*}
Suppose that $\overline{\varphi }(\mathcal{X}_{1}+\mathcal{X}_{2})\neq \overline{\varphi
}(\mathcal{X}_{1})+\overline{\varphi }(\mathcal{X}_{2}).$ In this case this means that $\overline{\varphi
}(\mathcal{X}_{1})+\overline{ \varphi }(\mathcal{X}_{2})\notin Im \overline{\varphi }.$ But by construction
$\overline{\overline{\varphi }(\mathcal{X}_{1}+\mathcal{X}_{2})} \in Im \overline{\varphi }$ and this coincides with
$\overline{\varphi }( \mathcal{X}_{1}+\mathcal{X}_{2}).$

\bigskip

\subsection{On the monotony of the product}

We define on $\mathcal{A}_{4}$ a partial order relation by
\begin{equation*}
\left\{
\begin{array}{l}
(x_{1},x_{2},0,0)\leq (y_{1},y_{2},0,0)\Longleftrightarrow y_{1}\leq x_{1} \text{ and }x_{2}\leq y_{2}, \\
(x_{1},x_{2},0,0)\leq (0,y_{2},y_{3},0)\Longleftrightarrow \text{ }x_{2}\leq y_{2}, \\ (0,x_{2},x_{3},0)\leq
(0,y_{2},y_{3},0)\Longleftrightarrow x_{3}\leq y_{3} \text{ and }x_{2}\leq y_{2}, \\ (0,0,x_{3},x_{4})\leq
(0,y_{2},y_{3},0)\Longleftrightarrow \text{ }x_{3}\leq y_{3}, \\ (0,0,x_{3},x_{4})\leq
(0,0,y_{3},y_{4})\Longleftrightarrow x_{3}\leq y_{3} \text{ and }y_{4}\leq x_{4}.
\end{array}
\right.
\end{equation*}

\begin{proposition}
\textbf{Monotony property: }Let $\mathcal{X}_{1},\mathcal{X}_{2}\in \overline{\mathbb{IR}}$. Then
\begin{equation*}
\left\{
\begin{array}{l}
\mathcal{X}_{1}\subset \mathcal{X}_{2}\Longrightarrow \mathcal{X}_{1}\bullet \mathcal{Z}\subset \mathcal{X}_{2}\bullet
\mathcal{Z}\text{ for all } \mathcal{Z}\in \overline{\mathbb{IR}}. \\ \overline{\varphi }(\mathcal{X}_{1})\leq
\overline{\varphi }(\mathcal{X} _{2})\Longrightarrow \overline{\varphi }(\mathcal{X}_{1}\bullet \mathcal{Z} )\leq
\overline{\varphi }(\mathcal{X}_{2}\bullet \mathcal{Z})
\end{array}
\right.
\end{equation*}
\end{proposition}

\noindent \textit{Proof. }Let us note that the second property is equivalent to the first.\ It is its translation in
$\overline{\mathcal{A}_{4}}.$ For any $\mathcal{X}\in \overline{\mathbb{IR}},$ we denote by $\overline{
\mathcal{X}}=\overline{\varphi }(\mathcal{X)}$ its image in $\overline{ \mathcal{A}_{4}}.$ Let
$\mathcal{X}_{1}$,$\mathcal{X}_{2}$ and $\mathcal{Z}$ be in $\overline{\mathbb{IR}}$. We denote by
$(z_{1},z_{2,}z_{3},z_{4})$ the image of $\mathcal{Z}$ in $\overline{\mathcal{A}_{4}},$ this implies $ z_{i}\geq 0$ and
$z_{1}z_{3}=z_{2}z_{4}=0.$ We assume that $\overline{ \varphi }(\mathcal{X}_{1})\leq \overline{\varphi
}(\mathcal{X}_{2}).$

\bigskip

\noindent \textit{First case. }$\overline{\mathcal{X}_{1}}=(x_{1},x_{2},0,0),
\overline{\mathcal{X}_{2}}=(y_{1},y_{2},0,0).$ We have
\begin{equation*}
\left\{
\begin{array}{l}
\overline{\varphi }(\mathcal{X}_{1}\bullet \mathcal{Z} )=(x_{1}z_{1},x_{2}z_{2},x_{2}z_{3},x_{1}z_{4}), \\
\overline{\varphi }(\mathcal{X}_{2}\bullet \mathcal{Z} )=(y_{1}z_{1},y_{2}z_{2},y_{2}z_{3},y_{1}z_{4}).
\end{array}
\right.
\end{equation*}
As $y_{1}z_{1}\leq x_{1}z_{1}$ and $x_{2}z_{2}\leq y_{2}z_{2},$ then $ \overline{\varphi }(\mathcal{X}_{1}\bullet
\mathcal{Z})\leq \overline{ \varphi }(\mathcal{X}_{2}\bullet \mathcal{Z}).$

\noindent \textit{Second case.} $\overline{\mathcal{X}_{1}}
=(x_{1},x_{2},0,0),\overline{\mathcal{X}_{2}}=(0,y_{2},y_{3},0).$ We have
\begin{equation*}
\left\{
\begin{array}{l}
\overline{\varphi }(\mathcal{X}_{1}\bullet \mathcal{Z} )=(x_{1}z_{1},x_{2}z_{2},x_{2}z_{3},x_{1}z_{4}), \\
\overline{\varphi }(\mathcal{X}_{2}\bullet \mathcal{Z} )=(0,y_{2}z_{2}+y_{3}z_{3},y_{2}z_{3}+y_{3}z_{2},0).
\end{array}
\right.
\end{equation*}
If $z_{1}\neq 0$ we have $y_{2}z_{2}+y_{3}z_{3}\geq x_{2}z_{2}$. If $z_{1}=0$ and $z_{2}\neq 0$, we have
$y_{2}z_{3}+y_{3}z_{2}\geq x_{2}z_{3}$ and $ y_{2}z_{2}+y_{3}z_{3}\geq x_{2}z_{2}$. If $z_{1}=z_{2}=0$, we have $
x_{2}z_{3}\leq y_{2}z_{3}$. This implies that in any case $\overline{\varphi }(\mathcal{X}_{1}\bullet \mathcal{Z})\leq
\overline{\varphi }(\mathcal{X} _{2}\bullet \mathcal{Z}).$

\noindent \textit{Third case.} $\overline{\mathcal{X}_{1}}=(0,x_{2},x_{3},0),
\overline{\mathcal{X}_{2}}=(0,y_{2},y_{3},0)$. We have
\begin{equation*}
\left\{
\begin{array}{l}
\overline{\varphi }(\mathcal{X}_{1}\bullet \mathcal{Z} )=(0,x_{2}z_{2}+x_{3}z_{3},x_{2}z_{3}+x_{3}z_{2},0), \\
\overline{\varphi }(\mathcal{X}_{2}\bullet \mathcal{Z} )=(0,y_{2}z_{2}+y_{3}z_{3},y_{2}z_{3}+y_{3}z_{2},0).
\end{array}
\right.
\end{equation*}
Thus
\begin{equation*}
\overline{\varphi }(\mathcal{X}_{1}\bullet \mathcal{Z})\leq \overline{ \varphi }(\mathcal{X}_{2}\bullet
\mathcal{Z})\Longleftrightarrow \left\{
\begin{array}{c}
(x_{2}-y_{2})z_{2}+(x_{3}-y_{3})z_{3}\leq 0, \\ (x_{2}-y_{2})z_{3}+(x_{3}-y_{3})z_{2}\leq 0.
\end{array}
\right.
\end{equation*}
But $(x_{2}-y_{2})$, $(x_{3}-y_{3})\leq 0$ and $\ z_{2},z_{3}\geq 0$. This implies $\overline{\varphi
}(\mathcal{X}_{1}\bullet \mathcal{Z})\leq \overline{\varphi }(\mathcal{X}_{2}\bullet \mathcal{Z}).$

\bigskip \noindent \textit{Fourth case.} $\overline{\mathcal{X}_{1}}
=(0,0,x_{3},x_{4}),\overline{\mathcal{X}_{2}}=(0,y_{2},y_{3},0).$ We have
\begin{equation*}
\left\{
\begin{array}{l}
\overline{\varphi }(\mathcal{X}_{1}\bullet \mathcal{Z} )=(x_{4}z_{4},x_{3}z_{3},x_{3}z_{2},x_{4}z_{1}), \\
\overline{\varphi }(\mathcal{X}_{2}\bullet \mathcal{Z} )=(0,y_{2}z_{2}+y_{3}z_{3},y_{2}z_{3}+y_{3}z_{2},0).
\end{array}
\right.
\end{equation*}
If $z_{4}\neq 0$, then $z_{2}=0$ and $x_{3}z_{3}\leq y_{3}z_{3}$. If $ z_{4}=0 $ and $z_{3}\neq 0$, we have
$x_{3}z_{3}\leq y_{2}z_{2}+y_{3}z_{3}$ and $x_{3}z_{2}\leq y_{2}z_{3}+y_{3}z_{2}$. If $z_{4}=z_{3}=0$ then $
x_{3}z_{2}\leq x_{3}z_{2}$. This implies that in any case $\overline{\varphi }(\mathcal{X}_{1}\bullet \mathcal{Z})\leq
\overline{\varphi }(\mathcal{X} _{2}\bullet \mathcal{Z}).$

\noindent \textit{Fifth case.} $\overline{\mathcal{X}_{1}}=(0,0,x_{3},x_{4}),
\overline{\mathcal{X}_{2}}=(0,0,y_{3},y_{4})$. We have
\begin{equation*}
\left\{
\begin{array}{l}
\overline{\varphi }(\mathcal{X}_{1}\bullet \mathcal{Z} )=(x_{4}z_{4},x_{3}z_{3},x_{3}z_{2},x_{4}z_{1}), \\
\overline{\varphi }(\mathcal{X}_{2}\bullet \mathcal{Z} )=(y_{4}z_{4},y_{3}z_{3},y_{3}z_{2},y_{4}z_{1}).
\end{array}
\right.
\end{equation*}
But $\overline{\varphi }(\mathcal{X}_{1}\bullet \mathcal{Z})\leq \overline{ \varphi }(\mathcal{X}_{2}\bullet
\mathcal{Z})$ is equivalent to $\left\{
\begin{array}{c}
x_{3}z_{2}\leq y_{3}z_{2}, \\ y_{4}z_{1}\leq x_{4}z_{1},
\end{array}
\right. $ and this is always satisfied.

\subsection{\protect\bigskip Application}

To work in the intervals set, we propose the following program:

1. Translate the problem in $\mathcal{A}_{4}$ through $\overline{\mathbb{IR} }.$

2. Solve the problem in $\mathcal{A}_{4}$ (which is a normed associative algebra).

3. Return to $\overline{\mathbb{IR}}$ and then to $\mathbb{IR}$.

\medskip The last condition require the use of $\psi .$

\subsection{Remark}

In the first version of this paper we have considered an another product in $ \mathcal{A}_{4}$ but this product was not
minimal. \qquad

\section{Divisibility and an Euclidean division}

We have computed the invertible elements of $\mathcal{A}_{4}.$ If $x=( x_{1},x_{2},x_{3},x_{4})\in \mathcal{A}_{4}$ and
if $\Delta =(x_{1}^{2}-x_{4}^{2})(x_{2}^{2}-x_{3}^{2})\neq 0$ then
\begin{equation*}
x^{-1}=\left( \frac{x_{1}}{x_{1}^{2}-x_{4}^{2}},\frac{x_{2}}{
x_{2}^{2}-x_{3}^{2}},\frac{x_{3}}{x_{2}^{2}-x_{3}^{2}},\frac{x_{4}}{ x_{1}^{2}-x_{4}^{2}}\right) .
\end{equation*}
The elements associated to $\mathcal{X}=\overline{(K,0)}$ are of the form
\begin{equation*}
\left\{
\begin{array}{l}
(x_{1},x_{2},0,0)\text{ if }0<x_{1}<x_{2}, \\ (0,x_{2},-x_{1},0)\text{ if }x_{1}<0<x_{2}, \\ (0,0,-x_{1},-x_{2})\text{
if }x_{1}<x_{2}<0,
\end{array}
\right.
\end{equation*}
and to $\mathcal{X}\in \overline{(0,K)}$
\begin{equation*}
\left\{
\begin{array}{l}
(0,0,x_{1},x_{2})\text{ if }0<x_{1}<x_{2}, \\ (-x_{1},0,0,x_{2})\text{ if }x_{1}<0<x_{2}, \\ (-x_{1},-x_{2},0,0)\text{
if }x_{1}<x_{2}<0.
\end{array}
\right.
\end{equation*}
The inverse of $(x_1,x_2,0,0)$ with $0<x_{1}<x_{2}$ is $\left(\dfrac{1}{x_{1} },\dfrac{1}{x_{2}},0,0\right)$.

\noindent The inverse of $(0,0,-x_{1},-x_{2})\text{ with }x_{1}<x_{2}<0$ is $
\left(0,0,-\dfrac{1}{x_{1}},-\dfrac{1}{x_{2}}\right)$.

\noindent The inverse of $(0,0,x_{1},x_{2})\text{ with }0<x_{1}<x_{2}$ is $
\left(0,0,\dfrac{1}{x_{1}},\dfrac{1}{x_{2}}\right)$.

\noindent The inverse of $(-x_{1},-x_{2},0,0)\text{ with }x_{1}<x_{2}<0$ is $
\left(-\dfrac{1}{x_{1}},-\dfrac{1}{x_{2}},0,0\right)$.

\noindent For $\mathcal{X}=(0,x_2,-x_1,0)$ or $(-x_1,0,0,x_2)$ with $ x_1x_2<0 $, then $\Delta=0$ and $\mathcal{X}$ is
not invertible. Then if $ \Delta \neq 0$ the inverse is always represented by an element of $\overline{ \mathbb{IR}} $
throught $\psi .$

\subsection{\protect\bigskip Division by an invertible element}

We denote by $\overline{\mathbb{IR}}^{+}$ the subset $\overline{(X,0)}$ with $X=[x_{1},x_{2}]$ and $0\leq x_{1}.$

\begin{proposition}
Let $\mathcal{X}=\overline{(X,0)}$ and $\mathcal{Y}=\overline{(Y,0)}$ be in $ \overline{\mathbb{IR}}^{+}$ with
$X=[x_{1},x_{2}],$ $Y=[y_{1},y_{2}].$ If $\ \dfrac{y_{2}}{y_{1}}\geq \dfrac{x_{2}}{x_{1}}$ then there exists an unique
$ \mathcal{Z}=\overline{(Z,0})\in \overline{\mathbb{IR}}^{+}$ such that $ \mathcal{Y=XZ}.$
\end{proposition}

\noindent\textit{Proof}. Let $\mathcal{Z}$ be defined by $c(\mathcal{Z})=
\frac{1}{2}\left(\dfrac{y_{2}}{x_{2}}+\dfrac{y_{1}}{x_{1}}\right)$ and $l(
\mathcal{Z})=\left(\dfrac{y_{2}}{x_{2}}-\dfrac{y_{1}}{x_{1}}\right).$ Then $ l(\mathcal{Z})\geq 0$ if and only if
$\dfrac{y_{2}}{x_{2}}\geq \dfrac{y_{1}}{ x_{1}}$ that is $\dfrac{y_{2}}{y_{1}}\geq \dfrac{x_{2}}{x_{1}}.$ Thus we have
$\mathcal{Y=XZ}.$ In fact
\begin{equation*}
\left( \overline{\varphi }(\mathcal{X})\right)^{-1}=\overline{\left(\dfrac{1 }{x_{1}},\dfrac{1}{x_{2}},0,0\right)}=\psi
\left( \ \overline{(0,[-\dfrac{1}{ x_{1}},-\dfrac{1}{x_{2}}])}\ \right).
\end{equation*}
Thus
\begin{equation*}
\overline{\varphi }(\mathcal{Y})\cdot \overline{\varphi }(\mathcal{X} )^{-1}=(y_{1},y_{2},0,0)\cdot
\left(\dfrac{1}{x_{1}},\dfrac{1}{x_{2}} ,0,0\right)=\left(\dfrac{y_{1}}{x_{1}},\dfrac{y_{2}}{x_{2}},0,0\right).
\end{equation*}
As $\dfrac{y_{1}}{x_{1}}\leq \dfrac{y_{2}}{x_{2}},$
\begin{equation*}
\psi (\overline{\varphi }(\mathcal{Y})\cdot \overline{\varphi }(\mathcal{X}
)^{-1})=\overline{\left([\dfrac{y_{1}}{x_{1}},\dfrac{y_{2}}{x_{2}}],0\right)} .
\end{equation*}
We can note also that
\begin{equation*}
\overline{\left(0,[-\dfrac{1}{x_{1}},-\dfrac{1}{x_{2}}]\right)}\bullet
\overline{([y_{1},y_{2}],0)}=\overline{\left([\dfrac{y_{1}}{x_{1}},\dfrac{ y_{2}}{x_{2}}],0\right)}.
\end{equation*}
Then the divisibility corresponds to the multiplication by the inverse element.

\bigskip

\subsection{Division by a non invertible element}

Let $X=[-x_{1},x_{2}]$ with $x_{1},x_{2}>0.$ We have seen that $\varphi (X)=(0,x_{2},x_{1},0)$ is not invertible in
$\mathcal{A}_{4}.$ For any $ M=(y_{1},y_{2},y_{3},y_{4})\in \mathcal{A}_{4}$ we have
\begin{equation*}
\varphi (X)\cdot M=(0,x_{2}y_{2}+x_{1}y_{3},x_{1}y_{2}+x_{2}y_{3},0)
\end{equation*}
and this point represents a non invertible interval. Thus we can solve the equation $\mathcal{Y=X\bullet Z}$ for
$\mathcal{X}=\overline{ ([-x_{1},x_{2}],0)}$ , $\mathcal{Y}=\overline{([-y_{1},y_{2}],0)}$ with $ x_{1},x_{2}>0$ and
$y_{1},y_{2}>0.$ Putting $\varphi (\mathcal{Z} )=(z_{1},z_{2},z_{3},z_{4})$, we obtain
\begin{equation*}
(0,y_{2},y_{1},0)=(0,x_{2},x_{1},0)\cdot (z_{1},z_{2},z_{3},z_{4}),
\end{equation*}
that is
\begin{equation*}
\left\{
\begin{array}{c}
y_{2}=x_{2}z_{2}+x_{1}z_{3}, \\ y_{1}=x_{2}z_{3}+x_{1}z_{2},
\end{array}
\right.
\end{equation*}
or
\begin{equation*}
\left(
\begin{array}{c}
y_{1} \\ y_{2}
\end{array}
\right) =\left(
\begin{array}{cc}
x_{1} & x_{2} \\ x_{2} & x_{1}
\end{array}
\right) \left(
\begin{array}{c}
z_{2} \\ z_{3}
\end{array}
\right) .
\end{equation*}
If $x_{1}^{2}-x_{2}^{2}\neq 0,$
\begin{equation*}
\left\{
\begin{array}{l}
z_{2}=\dfrac{x_{1}y_{1}-x_{2}y_{2}}{x_{1}^{2}-x_{2}^{2}}, \\
z_{3}=\dfrac{-x_{2}y_{1}+x_{1}y_{2}}{x_{1}^{2}-x_{2}^{2}}.
\end{array}
\right.
\end{equation*}
If $x_{1}^{2}-x_{2}^{2}=0$ then $x_{1}=x_{2}$ and the center of $ X=[-x_{1},x_{1}]$ is $0.$ Let us assume that
$x_{1}\neq x_{2}.$ If $ x_{1}^{2}-x_{2}^{2}<0$ $\ $that is $x_{1}<x_{2}$ then
\begin{equation*}
\left\{
\begin{array}{c}
x_{1}y_{1}-x_{2}y_{2}<0, \\ x_{1}y_{2}-x_{2}y_{1}<0,
\end{array}
\right.
\end{equation*}
and $\dfrac{x_{1}}{x_{2}}<\dfrac{y_{2}}{y_{1}},$ $\dfrac{x_{1}}{x_{2}}< \dfrac{y_{1}}{y_{2}}.$ If $\alpha
=\dfrac{x_{1}}{x_{2}}<1$ we have $ y_{2}>\alpha y_{1},$ $y_{1}>\alpha y_{2}$ then $y_{2}>\alpha ^{2}y_{2}$ and $
1-\alpha ^{2}>0.$ This case admits solution.

\begin{proposition}
Let $\mathcal{X}=\overline{([-x_{1},x_{2}],0)}$ with $x_{1},x_{2}>0$ and $ x_{1}<x_{2}.$ Then for any
$\mathcal{Y}=\overline{([-y_{1},y_{2}],0)}$ with $ y_{1},y_{2}>0$ and $\dfrac{x_{1}}{x_{2}}<\dfrac{y_{2}}{y_{1}},$
$\dfrac{x_{1} }{x_{2}}<\dfrac{y_{1}}{y_{2}},$ there is $\mathcal{Z}=\overline{ ([-z_{1},z_{2}],0)}$ such that
$\mathcal{Y=X\bullet Z}$.
\end{proposition}

\medskip

\noindent Suppose now that $x_{1}^{2}-x_{2}^{2}>0$ that us $x_{1}>x_{2}.$ In this case we have
\begin{equation*}
\left\{
\begin{array}{c}
x_{1}y_{1}-x_{2}y_{2}>0, \\ x_{1}y_{2}-x_{2}y_{1}>0,
\end{array}
\right.
\end{equation*}
that is $\dfrac{y_{2}}{y_{1}}<\dfrac{x_{1}}{x_{2}}$ and $\dfrac{y_{1}}{y_{2}} <\dfrac{x_{1}}{x_{2}}.$

\begin{proposition}
Let $\mathcal{X}=\overline{([-x_{1},x_{2}],0)}$ with $x_{1},x_{2}>0$ and $ x_{1}>x_{2}.$ For any
$\mathcal{Y}=\overline{([-y_{1},y_{2}],0)}$ with $ y_{1},y_{2}>0,$ $\dfrac{x_{1}}{x_{2}}>\dfrac{y_{2}}{y_{1}},$
$\dfrac{x_{1}}{ x_{2}}>\dfrac{y_{1}}{y_{2}},$ there is $\mathcal{Z}=\overline{ ([-z_{1},z_{2}],0)}$ such that
$\mathcal{Y=X\bullet Z}$.
\end{proposition}

\noindent \textbf{Example.} $\mathcal{X}=\overline{([-4,2],0)},$ $\mathcal{Y} =\overline{([-2,3],0)}.$ We have
$\dfrac{x_{2}}{x_{1}}=\dfrac{1}{2},$ $ \dfrac{x_{1}}{x_{2}}=2$ and $\dfrac{3}{2}<2<6.$ Then $\mathcal{Z}$ exists and it
is equal to $\mathcal{Z}=\overline{([-\dfrac{8}{12},\dfrac{2}{12}],0)} .$

\subsection{\protect\bigskip An Euclidean division}

Consider $\mathcal{X}=\overline{([x_{1},x_{2}],0)}$ and $\mathcal{Y}= \overline{([y_{1},y_{2}],0)\text{ }}$in
$\overline{\mathbb{IR}}^{+}.$ We have seen that $\mathcal{Y}$ is divisible by $\mathcal{X}$ as soon as $
\dfrac{x_{1}}{x_{2}}\geq \dfrac{y_{1}}{y_{2}}.$ We suppose now that $\dfrac{ x_{1}}{x_{2}}<\dfrac{y_{1}}{y_{2}}.$ In
this case we have

\begin{theorem}
Let $\mathcal{X}$ and $\mathcal{Y}$ be in $\overline{\mathbb{IR}}^{+}$ with $
\dfrac{x_{1}}{x_{2}}<\dfrac{y_{1}}{y_{2}}.$ There is a unique pair $( \mathcal{Z}$,$\mathcal{R})$ unique in
$\overline{\mathbb{IR}}^{+}$ such that
\begin{equation*}
\left\{
\begin{array}{l}
\mathcal{Y=X\bullet Z}+\mathcal{R}, \\ l(\mathcal{R})=0\text{ and }c(\mathcal{R)}\text{ minimal.}
\end{array}
\right.
\end{equation*}
This pair is given by
\begin{equation*}
\left\{
\begin{array}{l}
\mathcal{Z}=\dfrac{y_{2}-y_{1}}{x_{2}-x_{1}} \ \overline{([1,1],0)}, \\
\mathcal{R}=\dfrac{x_{2}y_{1}-x_{1}y_{2}}{x_{2}-x_{1}}\ \overline{([1,1],0)}.
\end{array}
\right.
\end{equation*}
\end{theorem}

\noindent \textit{Proof.} We consider $\mathcal{Z=}$ $\overline{ ([z_{1},z_{2}],0)}$ with $z_{1}>0.$ Then
$\mathcal{Y=X\bullet Z}+\mathcal{R}$ gives
\begin{equation*}
\mathcal{R}=\overline{([y_{1},y_{2}],[z_{1}x_{1},z_{2}x_{2}])}.
\end{equation*}
We have $\mathcal{R}\in \overline{\mathbb{IR}}^{+}$ if and only if \ $0\leq y_{1}-x_{1}z_{1}\leq y_{2}-x_{2}z_{2}$ that
is
\begin{equation*}
\left\{
\begin{array}{l}
z_{1}\leq \dfrac{y_{1}}{x_{1}}, \\ z_{2}\leq \dfrac{y_{2}}{x_{2}}, \\ z_{1}\geq \dfrac{y_{1}-y_{2}+x_{2}z_{2}}{x_{1}}.
\end{array}
\right.
\end{equation*}
The condition $z_{1}\leq z_{2}$ implies $\dfrac{y_{1}-y_{2}+x_{2}z_{2}}{x_{1} }\leq z_{2}$ that is $z_{2}\leq
\dfrac{y_{2}-y_{1}}{x_{2}-x_{1}}.$ Consider the case $z_{2}=\dfrac{y_{2}-y_{1}}{x_{2}-x_{1}}.$ Then $z_{1}\geq \dfrac{
y_{1}-y_{2}+x_{2}z_{2}}{x_{1}}=\dfrac{y_{2}-y_{1}}{x_{2}-x_{1}}=z_{2}$ and $ z_{1}=z_{2}.$ This case corresponds to
\begin{equation*}
\left\{
\begin{array}{l}
\mathcal{Z}=\dfrac{y_{2}-y_{1}}{x_{2}-x_{1}}\ \overline{([1,1],0)}, \\
\mathcal{R}=\dfrac{x_{2}y_{1}-x_{1}y_{2}}{x_{2}-x_{1}}\ \overline{([1,1],0)}.
\end{array}
\right.
\end{equation*}
Let us note that $y_{1}x_{2}-x_{1}y_{2}>0$ is equivalent to $\dfrac{y_{1}}{ y_{2}}>\dfrac{x_{1}}{x_{2}}$ which is
satisfied by hypothesis. We have also for this solution $l(\mathcal{R})=0$ and $c(\mathcal{R})=\dfrac{
x_{2}y_{1}-x_{1}y_{2}}{x_{2}-x_{1}}.$

\noindent Conversely, if $l(\mathcal{R})=0$ then $ y_{1}-z_{1}x_{1}=y_{2}-z_{2}x_{2}$ and
$z_{1}=z_{2}\dfrac{x_{2}}{x_{1}}+ \dfrac{y_{1}-y_{2}}{x_{1}}.$ As $z_{1}>0 $, we obtain $z_{2}>\dfrac{
y_{1}-y_{2}}{x_{1}}$ and $z_{1}\leq z_{2}$ implies
\begin{equation*}
\dfrac{y_{2}-y_{1}}{x_{2}}\leq z_{2}\leq \dfrac{y_{2}-y_{1}}{x_{2}-x_{1}}.
\end{equation*}
But $c(\mathcal{R})=y_{1}-z_{1}x_{1}=y_{2}-z_{2}x_{2}.$ Thus
\begin{equation*}
\dfrac{x_{2}y_{1}-x_{1}y_{2}}{x_{2}-x_{1}}\leq c(\mathcal{R})\leq y_{1}.
\end{equation*}
The norm is minimal when $c(\mathcal{R})=\dfrac{x_{2}y_{1}-x_{1}y_{2}}{ x_{2}-x_{1}}.$ \

\noindent\textbf{Example.} Let $\mathcal{X=}(\overline{[1,4],0})$ and $ \mathcal{Y}=(\overline{[1,3],0}).$ We have
$\dfrac{x_{1}}{x_{2}}=\dfrac{1}{4} <\dfrac{y_{1}}{y_{2}}=\dfrac{1}{3}.$ Thus $\mathcal{Z}=\dfrac{2}{3}(
\overline{[1,1],0})$ and $\mathcal{R}=\dfrac{1}{3}(\overline{[1,1],0}).$ The division writes
\begin{equation*}
(\overline{[1,3],0})=(\overline{[1,4],0})\cdot (\overline{[\dfrac{2}{3},
\dfrac{2}{3}],0})+(\overline{[\dfrac{1}{3},\dfrac{1}{3}],0}).
\end{equation*}

\bigskip Suppose now that $\mathcal{X}$ and $\mathcal{Y}$ are not
invertible, that is $\mathcal{X}=\overline{([-x_{1},x_{2}],0)}$ and $ \mathcal{Y}=\overline{([-y_{1},y_{2}],0)}$ with
$x_{1,}x_{2},y_{1},y_{2}$ positive$.$ We have seen that $\mathcal{Y}$ is divisible by $\mathcal{X}$ as soon as
\begin{equation*}
\left\{
\begin{array}{l}
\dfrac{x_{1}}{x_{2}}>\dfrac{y_{2}}{y_{1}}\text{ and }\dfrac{x_{1}}{x_{2}}> \dfrac{y_{1}}{y_{2}}, \\ \text{or } \\
\dfrac{x_{1}}{x_{2}}<\dfrac{y_{2}}{y_{1}}\text{ and }\dfrac{x_{1}}{x_{2}}< \dfrac{y_{1}}{y_{2}}.
\end{array}
\right.
\end{equation*}
We suppose now that these conditions are not satisfied. For example we assume that
\begin{equation*}
\dfrac{x_{1}}{x_{2}}>\dfrac{y_{2}}{y_{1}}\text{ and }\dfrac{x_{1}}{x_{2}}< \dfrac{y_{1}}{y_{2}}
\end{equation*}
(The other case is similar). If $\mathcal{Y=X\bullet Z}+\mathcal{R}$ then $ \mathcal{R}=(\overline{[-r_{1},r_{2}],0})$
with $r_{1}\geq 0$ and with $ r_{2}\geq 0$ because $\overline{\varphi }(\mathcal{R})=(0,r_{2},r_{1},0).$ This shows
that we can choose $\mathcal{Z}$ such that $\overline{\varphi }( \mathcal{R})=(0,z_{2},z_{3},0)$ and
\begin{equation*}
\left\{
\begin{array}{c}
z_{2}=\dfrac{x_{1}(y_{1}-r_{1})-x_{2}(y_{2}-r_{2})}{x_{1}^{2}-x_{2}^{2}}, \\
z_{3}=\dfrac{x_{1}(y_{2}-r_{2})-x_{2}(y_{1}-r_{1})}{x_{1}^{2}-x_{2}^{2}},
\end{array}
\right.
\end{equation*}
with the condition $z_{2}\geq 0$ and $z_{3}\geq 0$. If $x_{1}<x_{2}$ then this is equivalent to
\begin{equation*}
\left\{
\begin{array}{l}
\dfrac{x_{1}}{x_{2}}<\dfrac{y_{2}-r_{2}}{y_{1}-r_{1}}, \\ \dfrac{x_{1}}{x_{2}}<\dfrac{y_{1}-r_{1}}{y_{2}-r_{2}}.
\end{array}
\right.
\end{equation*}
If we suppose that $\mathcal{R\leq Y}$, thus $0<r_{2}<y_{2}$ and $ 0<r_{1}<y_{1},$ we obtain
\begin{equation*}
r_{1}>\dfrac{x_{2}}{x_{1}}r_{2}+\dfrac{-x_{2}y_{2}+x_{1}y_{1}}{x_{1}}<r_{1}<
\dfrac{x_{1}}{x_{2}}r_{2}+\dfrac{x_{2}y_{1}-x_{1}y_{2}}{x_{2}}.
\end{equation*}
Then lenght $l(\mathcal{R})=r_{1}+r_{2}$ is minimal if and only if $r_{2}=0$ and in this case
$r_{1}=\dfrac{x_{1}y_{1}-x_{2}y_{2}}{x_{1}}.$ We obtain
\begin{equation*}
\left\{
\begin{array}{l}
z_{2}=0, \\ z_{3}=\dfrac{y_{2}}{x_{1}}.
\end{array}
\right.
\end{equation*}

\begin{theorem}
Let $\mathcal{X}=\overline{([-x_{1},x_{2}],0)}$ with $x_{1},x_{2}>0$ and $ x_{1}>x_{2}.$ If
$\mathcal{Y}=\overline{([-y_{1},y_{2}],0)}$ with $ y_{1},y_{2}>0,$ satisfies
$\dfrac{x_{1}}{x_{2}}>\dfrac{y_{2}}{y_{1}}$ and $ \dfrac{x_{1}}{x_{2}}<\dfrac{y_{1}}{y_{2}},$ there is a unique pair
$\mathcal{ R} $,$\mathcal{Z}$ of non invertible elements such that
\begin{equation*}
\left\{
\begin{array}{l}
l(\mathcal{R})\text{ minimal,} \\ \mathcal{R<Y}\text{.}
\end{array}
\right.
\end{equation*}
This pair is given by
\begin{equation*}
\left\{
\begin{array}{l}
\mathcal{Z}=\overline{([-\dfrac{y_{2}}{x_{1}},0],0)}, \\
\mathcal{R}=\overline{([-\dfrac{x_{1}y_{1}-x_{2}y_{2}}{x_{1}},0],0)}.
\end{array}
\right.
\end{equation*}
\end{theorem}

\bigskip

\section{\protect\bigskip Applications}

\subsection{\protect\bigskip Differential calculus on $\overline{\mathbb{IR}}
$}

As $\overline{\mathbb{IR}}$ is a Banach space, we can describe a notion of differential function on it. \ Consider
$\mathcal{X} _{0}=\overline{(X_{0},0) }$ in $\overline{\mathbb{IR}}$ . The norm $||.||$ defines a topology on $
\overline{\mathbb{IR}}$ whose a basis of neighborhoods is given by the balls $\mathcal{B}(X_{0},\varepsilon )=\{X\in
\overline{\mathbb{IR}},|| \mathcal{X} \smallsetminus \mathcal{X} _{0}||<\varepsilon \}.$ Let us characterize the
elements of $\mathcal{B}(X_{0},\varepsilon ).$ $\mathcal{X} _{0}=\overline{(X_{0},0)}=\overline{([a,b],0)}.$

\begin{proposition}
Consider $\mathcal{X} _{0}=\overline{(X_{0},0)}=\overline{([a,b],0)}$ in $ \overline{\mathbb{IR}}$. Then every element
of $\mathcal{B} (X_{0},\varepsilon)$ is of type $\mathcal{X} =\overline{(X,0)}$ and satisfies
\begin{equation*}
l(X)\in B_{\mathbb{R}}(l(X_{0}),\varepsilon _{1})\text{ and }c(X)\in B_{ \mathbb{R}}(c(X_{0}),\varepsilon _{2})
\end{equation*}
with $\varepsilon _{1},\varepsilon_{2}\geq 0$ and $\varepsilon_{1}+\epsilon _{2}\leq \varepsilon,$ where
$B_{\mathbb{R}}(x,a)$ is the canonical open ball in $\mathbb{R}$\ of center $x$ and radius $a.$
\end{proposition}

\noindent\textit{Proof. }\textit{First case : }Assume that $\mathcal{X} = \overline{(X,0)}=\overline{([x,y],0)}$ . We
have
\begin{eqnarray*}
\mathcal{X} \smallsetminus \mathcal{X} _{0} &=&\overline{(X,X_{0})}= \overline{([x,y],[a,b])} \\ &=&\left\{
\begin{array}{c}
\overline{([x-a,y-b],0)}\text{ if }l(X)\geq l(X_{0}) \\ \overline{(0,[a-x,b-y])}\text{ if }l(X)\leq l(X_{0})
\end{array}
\right.
\end{eqnarray*}
If $l(X)\geq l(X_{0})$ we have
\begin{eqnarray*}
||\mathcal{X} \smallsetminus \mathcal{X} _{0}|| &=&(y-b)-(x-a)+\left\vert \frac{y-b+x-a}{2}\right\vert \\
&=&l(X)-l(X_{0})+|c(X)-c(X_{0})|.
\end{eqnarray*}
As \ $l(X)-l(X_{0})\geq 0$ and $|c(X)-c(X_{0})|\geq 0,$ each one of this term if less than $\varepsilon.$ If $l(X)\leq
l(X_{0})$ we have
\begin{equation*}
||\mathcal{X} \smallsetminus \mathcal{X} _{0}||=l(X_{0})-l(X)+|c(X_{0})-c(X)|.
\end{equation*}
and we have the same result.

\bigskip

\noindent\textit{Second case : }Consider $\mathcal{X} =\overline{(0,X)}= \overline{([x,y],0)}$ .\ We have
\begin{equation*}
\mathcal{X} \smallsetminus \mathcal{X} _{0}=\overline{(0,X_{0}+X)}=\overline{ ([x+a,y+b])}
\end{equation*}
and
\begin{equation*}
||\mathcal{X} \smallsetminus \mathcal{X} _{0}||=l(X_{0})+l(X)+|c(X_{0})+c(X)|.
\end{equation*}
In this case, we cannot have $||\mathcal{X} \smallsetminus \mathcal{X} _{0}||<\varepsilon$ thus $X\notin
\mathcal{B}(X_{0},\varepsilon).$

\medskip

\begin{definition}
A function $f:\overline{\mathbb{IR}}\longrightarrow \overline{\mathbb{R}}$ is continuous at $\mathcal{X} _{0}$ if
\begin{equation*}
\forall \varepsilon>0,\exists \eta >0\text{ such that }||\mathcal{X} \smallsetminus \mathcal{X} _{0}||<\varepsilon\text{
implies }||f(\mathcal{X} )\smallsetminus f(\mathcal{X} _{0})||<\varepsilon.
\end{equation*}
\end{definition}

Consider $(\mathcal{X} _{1},\mathcal{X} _{2})$ the basis of $\overline{ \mathbb{IR}}$ given in section 2. We have
\begin{equation*}
f(\mathcal{X} )=f_{1}(\mathcal{X} )\mathcal{X} _{1}+f_{2}(\mathcal{X} ) \mathcal{X} _{2}\text{ with
}f_{i}:\overline{\mathbb{IR}}\longrightarrow \mathbb{R}\text{.}
\end{equation*}
If $f$ is continuous at $\mathcal{X} _{0}$ so
\begin{equation*}
f(\mathcal{X} )\smallsetminus f(\mathcal{X} _{0})=(f_{1}(\mathcal{X} )-f_{1}( \mathcal{X} _{0}))\mathcal{X}
_{1}+(f_{2}(\mathcal{X} )-f_{2}(\mathcal{X} _{0}))\mathcal{X} _{2}.
\end{equation*}
To simplify notations let $\alpha =$ $f_{1}(\mathcal{X} )-f_{1}(\mathcal{X} _{0})$ and $\beta $ =$f_{2}(\mathcal{X}
)-f_{2}(\mathcal{X} _{0}).$ If $||f( \mathcal{X} )\smallsetminus f(\mathcal{X} _{0})||<\varepsilon$, and if we assume
$f_{1}(\mathcal{X} )-f_{1}(\mathcal{X} _{0})>0$ and $f_{2}(\mathcal{X} )-f_{2}(\mathcal{X} _{0})>0$ (other cases are
similar), then we have
\begin{equation*}
l(\alpha \mathcal{X} _{1}+\beta \mathcal{X} _{2})=l\overline{([\beta ,\alpha +\beta ],0)}<\varepsilon
\end{equation*}
thus $f_{1}(\mathcal{X} )-f_{1}(\mathcal{X} _{0})<\varepsilon.$ Simillary,
\begin{equation*}
c(\alpha \mathcal{X} _{1}+\beta \mathcal{X} _{2})=c\overline{([\beta ,\alpha +\beta ],0)}=\frac{\alpha }{2}+\beta
<\varepsilon
\end{equation*}
and this implies that $f_{2}(\mathcal{X} )-f_{2}(\mathcal{X} _{0})<\varepsilon.$

\begin{corollary}
$f$ is continuous at $\mathcal{X} _{0}$ if and only if $f_{1}$ and $f_{2}$ are continuous at $\mathcal{X} _{0}.\medskip
$
\end{corollary}

\noindent \textbf{Examples. }

\begin{itemize}
\item $f(\mathcal{X} )=\mathcal{X} .$ This function is continuous at any point.

\item $f(\mathcal{X} )=\mathcal{X} ^{2}.$ Consider $\mathcal{X} _{0}= \overline{(X_{0},0)}=\overline{([a,b],0)}$
    and $\mathcal{X} \in \mathcal{B} (X_{0},\varepsilon).$ We have
\begin{eqnarray*}
||\mathcal{X} ^{2}\smallsetminus \mathcal{X} _{0}^{2}|| &=&||(\mathcal{X} \smallsetminus \mathcal{X}
_{0})(\mathcal{X} +\mathcal{X} _{0})|| \\ &\leq &||\mathcal{X} \smallsetminus \mathcal{X} _{0}||||\mathcal{X} +
\mathcal{X} _{0}||.
\end{eqnarray*}
Given $\varepsilon>0,$ let $\eta =\dfrac{\varepsilon}{||\mathcal{X} + \mathcal{X} _{0}||},$ thus if $||\mathcal{X}
\smallsetminus \mathcal{X} _{0}||<\eta $, we have $||\mathcal{X} ^{2}\smallsetminus \mathcal{X}
_{0}^{2}||<\varepsilon$ and $f$ is continuous.

\item Consider $P=a_{0}+a_{1}X+\cdots +a_{n}X^{n}\in \mathbb{R[X]}$. We define
    $f:\overline{\mathbb{IR}}\longrightarrow \overline{\mathbb{IR}}$ with $f(\mathcal{X} )=a_{0}\mathcal{X}
    _{2}+a_{1}\mathcal{X} +\cdots +a_{n}^{n} \mathcal{X} ^{n}$ where $\mathcal{X} ^{n}=\mathcal{X} \cdot\mathcal{X}
    ^{n-1} $ . From the previous example, all monomials are continuous, it implies that $f$ is continuous.
\end{itemize}

\medskip

\begin{definition}
Consider $\mathcal{X}_{0}$ in $\overline{\mathbb{IR}}$ and $f:\overline{ \mathbb{IR}}\longrightarrow
\overline{\mathbb{IR}}$ continuous. We say that $ f$ is differentiable at $\mathcal{X}_{0}$ if there is
$g:\overline{\mathbb{IR }}\longrightarrow \overline{\mathbb{IR}}$ linear such as
\begin{equation*}
||f(\mathcal{X})\smallsetminus f(\mathcal{X}_{0})\smallsetminus g(\mathcal{X} \smallsetminus
\mathcal{X}_{0})||=o(||\mathcal{X}\smallsetminus \mathcal{X} _{0}||).
\end{equation*}
\end{definition}

\subsection{Study of the function $q_{2}$}

We consider the function $q_{2}:$ $\overline{\mathbb{IR}}\longrightarrow $ $ \overline{\mathbb{IR}}$ given by
\begin{equation*}
q_{2}\overline{([a,b],0)}=\left\{
\begin{array}{l}
\overline{([a^{2},b^{2}],0)}\text{ if }0<a<b, \\ \overline{([b^{2},a^{2}],0)}\text{ if }a<b<0, \\ \overline{([0,\sup
(a^{2},b^{2})])}\text{ if }a<0<b\text{.}
\end{array}
\right.
\end{equation*}
and $q_{2}\overline{(0,[a,b])}=q_{2}\overline{([a,b],0)}.$ For any invertible element $\mathcal{X}\in
\overline{\mathbb{IR}}$, we have $q_{2}( \mathcal{X)}=\mathcal{X\bullet X=X}^{2}.$ If $\mathcal{X}$ is not invertible,
it writes $\mathcal{X}=\overline{([a,b],0)}$ with $a<0<b$ ( we assume that $\mathcal{X}$ $\ $is of type
$\overline{(K,0)}$). In this case $ \mathcal{X}$ $\bullet \mathcal{X}$ $=\overline{([2ab,a^{2}+b^{2}],0)}$ and $
q_{2}\subset \mathcal{X}$ $\bullet \mathcal{X}$ .

\begin{proposition}
The function $q_{2}$ is continuous on $\overline{\mathbb{IR}}.$
\end{proposition}

\noindent \textit{Proof.} Let $\mathcal{X}_{0}\in \overline{\mathbb{IR}}.$ Assume that
$\mathcal{X}_{0}=\overline{([a,b],0)}$ with $0<a<b.$ An $\eta $ -neighborhood is represented by the parallelogram
$(A,B,C,D)$ with $A=(a- \dfrac{\eta }{2},b+\dfrac{\eta }{2}),$ $B=(a-\eta ,b-\eta ),$ $C=(a+\dfrac{ \eta
}{2},b-\dfrac{\eta }{2}),$ $D=(a+\eta ,b+\eta ).$ We have $q_{2}(
\mathcal{X}_{0})=\mathcal{X}_{0}^{2}=\overline{([a^{2},b^{2}],0)}.$ For any $ \varepsilon >0$ we consider the
$\varepsilon $-neighborhood of $q_{2}( \mathcal{X}_{0}).$ it is represented by the parallelogram $
(A_{1},B_{1},C_{1},D_{1})$ with $A_{1}=(a^{2}-\dfrac{\varepsilon }{2},b^{2}+ \dfrac{\varepsilon }{2}),$
$B_{1}=(a^{2}-\varepsilon ,b^{2}-\varepsilon ),$ $ C_{1}=(a^{2}+\dfrac{\varepsilon }{2},b^{2}-\dfrac{\varepsilon
}{2}),$ $ D_{1}=(a^{2}+\varepsilon ,b^{2}+\varepsilon ).$ If $\eta $ satisfy
\begin{equation*}
\left\{
\begin{array}{c}
2a\eta +\eta ^{2}<\dfrac{\varepsilon }{2}, \\ \eta ^{2}-2a\eta >-\dfrac{\varepsilon }{2}
\end{array}
\right.
\end{equation*}
the \ for every point of the $\eta $-neighborhood of $\mathcal{X}_{0}$, the image $q_{2}(\mathcal{X})$ is contained in
the $\varepsilon $ -neighborhood of $q_{2}(\mathcal{X}_{0}).$ If $a\neq 0$, as $\varepsilon $ is infinitesimal we have
$\eta =\dfrac{\varepsilon }{8a}.$ If $a=0,$ we have $\eta =\varepsilon .$ Then $q_{2}$ is continuous at the point
$\mathcal{X} _{0}.$ Is $\mathcal{X}_{0}=\overline{([a,b],0)}$ with $a<b<0 $, taking $\eta =-\dfrac{\varepsilon }{8a}$
we prove in a similar way the continuity at $ \mathcal{X}_{0}.$

\noindent Assume that $\mathcal{X}_{0}=\overline{([a,b],0)}$ with $a<0<b$ then
$q_{2}(\mathcal{X}_{0})=\overline{([0,\sup (a^{2},b^{2})])}.$ If $ \mathcal{X}=\overline{([x,y],0)}$ is an $\eta
$-neighborhood of $\mathcal{X }_{0}$ with $q_{2}(\mathcal{X})=\overline{([0,\sup ((x+\eta )^{2},(y+\eta )^{2})])}$ then
$a-\eta <x<a+\eta , b-\eta <y<b+\eta $ and we can find $\eta $ such that $\sup (a^{2},b^{2})-\dfrac{\varepsilon
}{2}<\sup ((x+\eta )^{2},(y+\eta )^{2})<\sup (a^{2},b^{2})+\dfrac{\varepsilon }{2}.$ Thus $q_{2} $ is also continuous in
this point. As $q_{2}\overline{(0,K)}=q_{2}\overline{ (K,0)}$, we have the continuity of any point.

\begin{theorem}
The function $q_2$ is not differentiable.
\end{theorem}
\textit{Proof. } The function $q_{2}$ is differentiable at the point $ \mathcal{X}_{0}$ if there is a linear map $L$
such that
\begin{equation*}
\displaystyle\lim_{||\mathcal{X}\smallsetminus \mathcal{X}_{0}||\rightarrow 0}\frac{||q_{2}(\mathcal{X})\smallsetminus
q_{2}(\mathcal{X} _{0})\smallsetminus L(\mathcal{X}\smallsetminus \mathcal{X}_{0})||}{|| \mathcal{X}\smallsetminus
\mathcal{X}_{0}||}=0.
\end{equation*}
We consider \noindent $L$ be the linear function given by
\begin{equation*}
L(\mathcal{X})=2\mathcal{X}_{0}\bullet (\mathcal{X}).
\end{equation*}
 We assume that $\mathcal{X}_{0}=(\overline{[a,b],0})$ with $0<a<b$. If $\mathcal{X}
$ is in an infinitesimal neighborhood of $\mathcal{X}_{0}$, then $\mathcal{X} =(\overline{[x,y],0})$ with $0<x<y$.

\begin{itemize}
\item If $0<x-a<y-b$
\begin{equation*}
\mathcal{X}\smallsetminus \mathcal{X}_{0}=(\overline{[x,y],[a,b]})=(\overline{[x-a,y-b],0})
\end{equation*}
Thus $L(\mathcal{X}\smallsetminus \mathcal{X}_{0})=2(\overline{[a,b],0})\bullet
(\overline{[x-a,y-b],0})=2(\overline{[a(x-a),b(y-b),0})$ and
\begin{equation*}
\begin{array}{ll}
q_{2}(\mathcal{X})\smallsetminus q_{2}(\mathcal{X}_{0})\smallsetminus L( \mathcal{X}\smallsetminus \mathcal{X}_{0})
& =(\overline{[x^{2},y^{2}],0}) \smallsetminus (\overline{[a^{2},b^{2}],0})\smallsetminus
2(\overline{[a(x-a),b(y-b),0}]), \\ & =(\overline{[x^{2}-a^{2},y^{2}-b^{2}],0})\smallsetminus
2(\overline{[a(x-a),b(y-b),0]}), \\ & =(\overline{[(x-a)^{2},(y-b)^{2}],0}).
\end{array}
\end{equation*}
We deduce
\begin{equation*}
\begin{array}{ll}
||q_{2}(\mathcal{X})\smallsetminus q_{2}(\mathcal{X}_{0})\smallsetminus L( \mathcal{X}\smallsetminus
\mathcal{X}_{0})|| & =(y-b)^{2}-(x-a)^{2}+| \displaystyle\frac{(y-b)^{2}+(x-a)^{2}}{2}\mid , \\ &
=\displaystyle\frac{3(y-b)^{2}-(x-a)^{2}}{2}. \\ &
\end{array}
\end{equation*}
Thus
\begin{equation*}
\frac{||q_{2}(\mathcal{X})\smallsetminus q_{2}(\mathcal{X} _{0})\smallsetminus L(\mathcal{X}\smallsetminus
\mathcal{X}_{0})||}{|| \mathcal{X}\smallsetminus \mathcal{X}_{0}||}=\frac{3(y-b)^{2}-(x-a)^{2}}{ 3(y-b)-(x-a)}.
\end{equation*}
Then  $\dfrac{||q_{2}(\mathcal{X})\smallsetminus q_{2}(\mathcal{X} _{0})\smallsetminus L(\mathcal{X}\smallsetminus
\mathcal{X}_{0})||}{|| \mathcal{X}\smallsetminus \mathcal{X}_{0}||}\leq \varepsilon $ is equivalent to
$$3(y-b)^{2}-(x-a)^{2}\leq \varepsilon (3(y-b)-(x-a)).$$ \ The following picture represent the $\varepsilon
$-neighborhood of $\mathcal{X}_{0}$ and the set defined by the previous inequalities (in case of
$\mathcal{X}_{0}=\overline{([1,2],0)}$ and $ \varepsilon=0,5$).

\noindent
\includegraphics{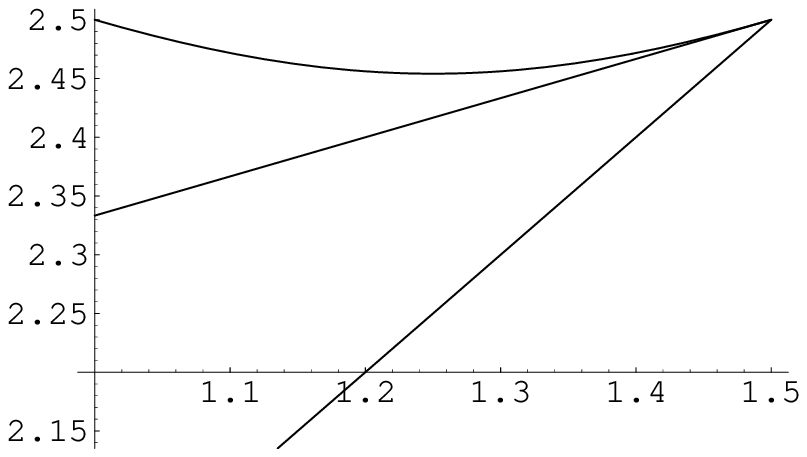}
\end{itemize}

\ Then for every point of the $\varepsilon $-neighborhood of $\mathcal{X} _{0}$, the $\varepsilon $ inequality of the
differentiability is satisfied. This shows that, if $q_2$ is differentiable at $\mathcal{X}_{0}$, then the differential
is the linear function $L(\mathcal{X})=2\mathcal{X}_{0} \bullet \mathcal{X}$.
\begin{itemize}
\item If $0<a-x<y-b.$ We find again the previous case.

\item If $0<y-b<a-x$, then
\begin{equation*}
\mathcal{X}\smallsetminus \mathcal{X}_{0}=(\overline{[x,y],[a,b]})=(\overline{[x-a,y-b],0}).
\end{equation*}
Thus $L(\mathcal{X}\smallsetminus \mathcal{X}_{0})=2(\overline{[a,b],0})\bullet
(\overline{[x-a,y-b],0})=2(\overline{[b(x-a),b(y-b),0})$ and
\begin{equation*}
\begin{array}{ll}
q_{2}(\mathcal{X})\smallsetminus q_{2}(\mathcal{X}_{0})\smallsetminus L( \mathcal{X}\smallsetminus \mathcal{X}_{0})
& =(\overline{[x^{2},y^{2}],0}) \smallsetminus (\overline{[a^{2},b^{2}],0})\smallsetminus
2(\overline{[b(x-a),b(y-b)],0}), \\ & =(\overline{[x^{2}-a^{2},y^{2}-b^{2}],0})\smallsetminus
2(\overline{[b(x-a),b(y-b)],0}), \\ & =(\overline{[(x-b)^{2}-(a-b)^{2},(y-b)^{2}],0}).
\end{array}
\end{equation*}
We deduce
\begin{equation*}
\frac{||q_{2}(\mathcal{X})\smallsetminus q_{2}(\mathcal{X} _{0})\smallsetminus L(\mathcal{X}\smallsetminus
\mathcal{X}_{0})||}{|| \mathcal{X}\smallsetminus \mathcal{X}_{0}||}=\frac{
3(y-b)^{2}-(x-b)^{2}+(a-b)^{2}}{(y-b)+3(a-x)}.
\end{equation*}
Then $\dfrac{||q_{2}(\mathcal{X})\smallsetminus q_{2}(\mathcal{X} _{0})\smallsetminus L(\mathcal{X}\smallsetminus
\mathcal{X}_{0})||}{|| \mathcal{X}\smallsetminus \mathcal{X}_{0}||}\leq \varepsilon $ is equivalent to
$$3(y-b)^{2}-(x-b)^{2}+(a-b)^{2}\leq \varepsilon (y-b)+3(a-x).$$ The following picture represent the $\varepsilon
$-neighborhood of $\mathcal{X} _{0}$ and the set $E$ defined  by the previous inequalities.

\noindent
\includegraphics{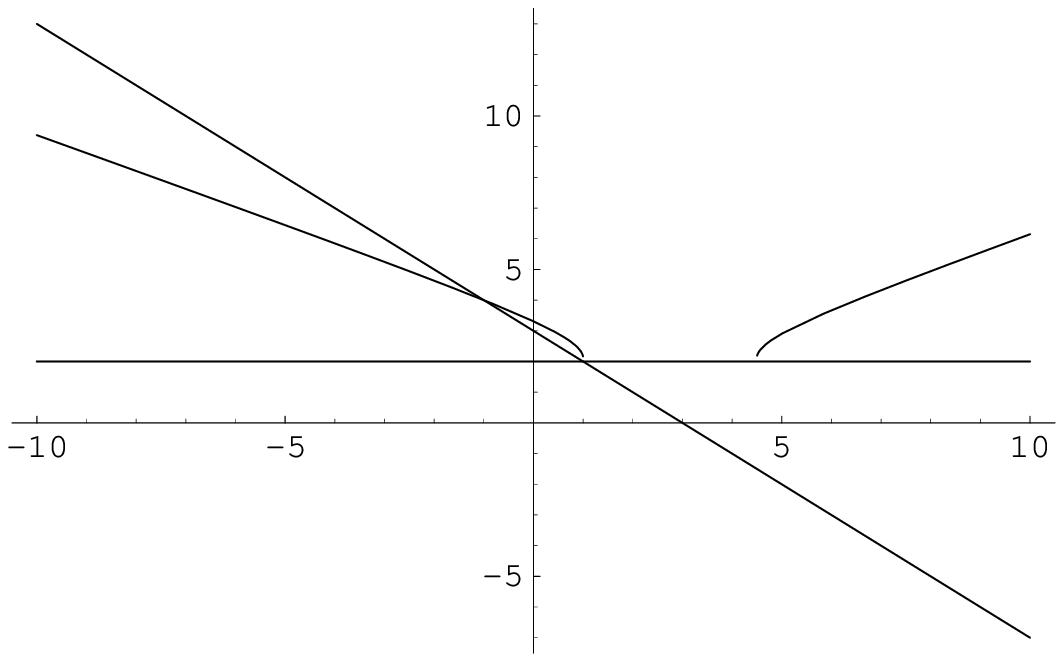}
\end{itemize}
We see that the representation of $E$  doesn't contains any points
 of the representation of a $\eta $-neighborhood of $\mathcal{X}_{0}$ for all $\eta$.
This gives a contradiction of the differentiability at  $\mathcal{X}_{0}$.

\subsection{Linear programming on $\overline{\mathbb{IR}}$}

We consider the following linear programm
\begin{equation*}
\left\{
\begin{array}{l}
Ax=b, \\ max (f(x)),
\end{array}
\right.
\end{equation*}
where $A$ is a real $n\times p$ matrix and $f:\mathbb{R}^{n}$ $ \longrightarrow \mathbb{R}$ a linear map.\ We want to
translate this problem in the context of intervals. This is easy because the set of intervals is endowed with a
vectorial space structure.\ Then we assume here that each variable $x_{i}$ belongs to an interval
$\mathcal{X}_{i}=\overline{(X_{i},0)} .$ We assume also that the constraints $b_{i}$ belongs to $Y_{i}=\overline{
(B_{i},0)}$ with $B_{i}=[b_{i}^{1},b_{i}^{2}]$ and $b_{i}^{1}\geq 0.$ We can extend the linear $f:\mathbb{R}^{n}$
$\longrightarrow \mathbb{R}$ which is written $f(x_{1},\cdots ,x_{n})=\sum a_{i}x_{i}$ by a linear map, denoted by
$\overline{f}$, given by $\overline{f}(\mathcal{X}_{1},\cdots ,\mathcal{X} _{n})=\sum a_{i}\mathcal{X}_{i}.$ In
$\overline{\mathbb{IR}}$ we claim that
\begin{equation*}
\mathcal{X}\geq \mathcal{X}^{\prime }\Longleftrightarrow \mathcal{X} \smallsetminus \mathcal{X}^{\prime
}=\overline{(K,0)}.
\end{equation*}
We have the following linear programming in terms of intervals:
\begin{equation*}
\left\{
\begin{array}{l}
A\mathcal{X}=\mathcal{B}\text{ with }\mathcal{X}\in \overline{\mathbb{IR}} ^{n}\text{ and }\mathcal{B\in
}\overline{\mathbb{IR}}^{p}, \\ \mathcal{B=(B}_{1},\cdots ,\mathcal{B}_{p})\text{ and }\mathcal{B}_{i}\geq 0,
\\
\max (\overline{f}(\mathcal{X})).
\end{array}
\right.
\end{equation*}
To solve this program, we extend the classical simplex algorithm. The simplex method is based on an algorithm using the
Gauss elimination. We consider the linear programm on the intervals, with $\mathcal{X}=(\mathcal{X} _{1},\cdots
,\mathcal{X}_{n})\in $ $\overline{\mathbb{IR}}^{n}$ and $ \mathcal{B=(B}_{1},\cdots ,\mathcal{B}_{p})\in $
$\overline{\mathbb{IR}} ^{p}. $ By hypothesis $\mathcal{B}_{i}=\overline{(Y_{i},0)}$ is a positive vector. We have to
define the pivot. The column pivot is defined by the largest positive coefficient of the economic function
$\overline{f.}$ In choosing the line, the goal is to keep the second member of positive vectors of
$\overline{\mathbb{IR}}^{p}$. Let $k$ be the number of the column containing the pivot. We note $l(Y_{j})$ the length of
the interval $Y_{j}$ and let $i$ be such that
\begin{equation*}
\dfrac{l(Y_{i})}{a_{ik}}=\min_{j}\left\{ \dfrac{l(Y_{j})}{a_{jk}} ,a_{jk}>0\right\} .
\end{equation*}
We choose $a_{i}^{k}$ as pivot. The line $\ l_{j}$ is transformed into $\ a_{ik}l_{j}-a_{jk}l_{i}.$ In this case the
second member becomes $a_{ik} \mathcal{B}_{j}\smallsetminus a_{jk}\mathcal{B}_{i}.$ Let us compute this value. We
assume $a_{jk}>0$
\begin{equation*}
\begin{array}{ll}
a_{ik}\mathcal{B}_{j}\smallsetminus a_{jk}\mathcal{B}_{i} & =\overline{ (a_{ik}Y_{j},a_{jk}Y_{i})} \\ &
=\overline{(a_{ik}[Y_{j}^{1},Y_{j}^{2}],a_{jk}[Y_{i}^{1},Y_{i}^{2}])} \\ & =\overline{
([a_{ik}Y_{j}^{1},a_{ik}Y_{j}^{2}],[a_{jk}Y_{i}^{1},a_{jk}Y_{i}^{2}])} \\ & =\overline{
([a_{ik}Y_{j}^{1}-a_{jk}Y_{i}^{1},a_{ik}Y_{j}^{2}-a_{jk}Y_{i}^{2}],0)}.
\end{array}
\end{equation*}
This is well defined. Indeed
\begin{eqnarray*}
(a_{ik}Y_{j}^{1}-a_{jk}Y_{i}^{1}) &<&a_{ik}Y_{j}^{2}-a_{jk}Y_{i}^{2} \\ &\Longleftrightarrow
&a_{jk}(Y_{i}^{2}-Y_{i}^{1})<a_{ik}(Y_{j}^{2}-Y_{j}^{1})
\\
&\Longleftrightarrow &\frac{l(Y_{i})}{a_{ik}}<\dfrac{l(Y_{j})}{a_{jk}},
\end{eqnarray*}
what is assumed by hypothesis. So $a_{ik}\mathcal{B}_{j}\smallsetminus a_{jk} \mathcal{B}_{i}$ is positive. In the case
where $a_{jk}<0$, $l(Y_{j})$ is transformed into $a_{ik}l(Y_{j})-a_{jk}l(Y_{i})$ and thus the second member is given
by
\begin{equation*}
\overline{(a_{ik}Y_{j},a_{jk}Y_{i})}=\overline{(a_{ik}Y_{j},0)}+\overline{
(-a_{jk}Y_{i},0)}=\overline{(a_{ik}Y_{j}-a_{jk}Y_{i},0)}.
\end{equation*}
\noindent This vector is always positive. This process, as in the classical simplex algorithm, gives in terms of
intervals the maximum of the economic function.

\subsection{Non Linear programming}

We consider here the following programm
\begin{equation*}
\left\{
\begin{array}{l}
Ax=b, \\ max (f(x)),
\end{array}
\right.
\end{equation*}
where $A$ is a real $n\times p$ matrix and $f:\mathbb{R}^{n}$ $ \longrightarrow \mathbb{R}$ a non linear differentiable
map. We translate this programm in terms of intervals:
\begin{equation*}
\left\{
\begin{array}{l}
A\mathcal{X}=\mathcal{B}\text{ with }\mathcal{X}\in \overline{\mathbb{IR}} ^{n}\text{ and }\mathcal{B\in
}\overline{\mathbb{IR}}^{p}, \\ \mathcal{B=(B}_{1},\cdots ,\mathcal{B}_{p})\text{ and }\mathcal{B}_{i}\geq 0,
\\
\max (\widetilde{f}(\mathcal{X})).
\end{array}
\right.
\end{equation*}
where $\widetilde{f}:\overline{\mathbb{IR}}^{n}$ $\longrightarrow \overline{ \mathbb{IR}}$ is the transferred
function.\ It is differentiable in the previous sense.\ Assume that the system $A\mathcal{X}=\mathcal{B}$ is
indeterminate and all the solutions are expressed in term of one parameter $ \mathcal{X}_{0}.$ In this case, a solution
of the non linear programming in given by a root of the equation $f^{\prime }(\mathcal{X}_{0})=0.$

\end{document}